\documentclass{gtart}

\def\ifplaintex{\expandafter\ifx\csname documentclass\endcsname\relax}

\ifplaintex 
\else
\fi

\expandafter\ifx\csname beginpicture\endcsname\relax
\expandafter\ifx\csname documentclass\endcsname\relax
\input pictex \else
\input prepictex \input pictex \input postpictex \fi\fi

\def\gt{{\mathsurround=0pt\it $\cal G\mskip-2mu$eometry \&\ 
$\cal T\!\!$opology}}        

\def\gtp{{\mathsurround=0pt\it $\cal G\mskip-2mu$eometry \&\ 
$\cal T\!\!$opology $\cal P\!$ublications}}  

\def\lognumber#1{\def\thelognumber{#1}}
\def\volumenumber#1{\def\thevolumenumber{#1}}
\def\papernumber#1{\def\thepapernumber{#1}}
\def\volumeyear#1{\def\thevolumeyear{#1}}

\def\pagenumbers#1#2{\def\startpage{#1}\def\finishpage{#2}}
\def\published#1{\def\publishdate{#1}}
\def\proposed#1{\def\theproposer{#1}}
\def\seconded#1{\def\theseconders{#1}}
\def\received#1{\def\receiveddate{#1}}
\def\revised#1{\def\reviseddate{#1}}
\def\accepted#1{\def\accepteddate{#1}}
\def\asciititle#1{\def\theasciititle{#1}}

\long\def\asciiabstract#1{\long\def\theasciiabstract{#1}}
\def\asciikeywords#1{\def\theasciikeywords{#1}}

\let\\\par\let\thelognumber\relax
\let\thevolumenumber\relax\let\thepapernumber\relax
\let\thevolumeyear\relax\let\thesamplenumber\relax\let\startpage\relax
\let\finishpage\relax\let\publishdate\relax\let\receiveddate\relax
\let\reviseddate\relax\let\accepteddate\relax\let\theasciititle\relax
\let\theasciiauthors\relax
\let\theasciiabstract\relax\let\theasciikeywords\relax
\let\theasciiemail\relax\let\theshortauthors\relax\let\theshorttitle\relax

\long\def\maketitlep{   

\count0=\startpage

\gt\hfill      
\beginpicture
\setcoordinatesystem units <0.33truein, 0.33truein> point at 2.2 0.9
\setplotsymbol ({$\cal G$})
\plotsymbolspacing=9truept
\circulararc 315 degrees from 0 1 center at 0 0
\setplotsymbol ({$\cal T$})
\circulararc 315 degrees from 1 -1 center at 1 0
\endpicture
\break
{\small\ifx\thesamplenumber\relax 
Volume \else Sample
\fi\thevolumenumber\ (\thevolumeyear)
\startpage--\finishpage\nl
Published: \publishdate}
\vglue 0.5truein plus 0.4fil minus 0.1truein

{\parskip=0pt\leftskip 0pt plus 1fil\def\\{\par\smallskip}{\ifplaintex\large
\else\Large\fi\bf\thetitle}\par\medskip}   

\vglue 0pt plus 0.1fil 

{\parskip=0pt\leftskip 0pt plus 1fil\def\\{\par}{\sc\theauthors}
\par\medskip}

\vglue 0pt plus 0.1fil 

{\small\parskip=0pt\let\newline\\
{\leftskip 0pt plus 1fil\def\\{\par}{\sl\theaddress}\par}
\expandafter\ifx\theemail\relax    
\relax\else\vglue 5pt plus 0.02fil minus 2pt\def\\{\stdspace{\rm 
and}\stdspace} 
\cl{Email:\stdspace\tt\theemail}\fi
\ifx\theurl\relax                  
\relax\else\vglue 5pt plus 0.02fil minus 2pt\def\\{\stdspace{\rm 
and}\stdspace}
\cl{URL:\stdspace\tt\theurl}\fi\par}

\vglue 7pt plus 0.3fil minus 3pt

{\bf Abstract}
\vglue 5pt plus 0.1fil minus 2pt

\theabstract

\vglue 7pt plus 0.3fil minus 3pt

{\bf AMS Classification numbers}\quad Primary:\quad \theprimaryclass

Secondary:\quad \thesecondaryclass

\vglue 5pt plus 0.3fil minus 2pt

{\bf Keywords}\quad \thekeywords

\vglue 10pt plus 0.5fil minus 5pt

{\small  Proposed: \theproposer\hfill Received: \receiveddate\nl
Seconded: \theseconders\hfill 
\ifx\reviseddate\relax                         
Accepted: \accepteddate                        
\else
Revised: \reviseddate                          
\fi}
\eject
}       

\let\maketitlepage\maketitlep
\let\maketitle\maketitlepage

\font\phead=cmsl9 scaled 950
\font=cmsl9 scaled 1050
\font\pnum=cmbx10 scaled 913
\font\lnum=cmbx10 
\font\pfoot=cmsl9 scaled 950
\font=cmsl9 scaled 1050
\ifplaintex
\headline{\vbox to 0pt{\vskip -4.5mm\line{\small\phead\ifnum
\count0=\startpage ISSN 1364-0380 (on line)
1465-3060 (printed) \hfill {\pnum\folio}\else\ifodd\count0\def\\{ }%
\ifx\theshorttitle\relax\thetitle\else\theshorttitle\fi\hfill{\pnum\folio}
\else\def\\{ and }{\pnum\folio}\hfill\ifx\theshortauthors\relax\theauthors
\else\theshortauthors\fi\fi\fi}\vss}}
\footline{\vbox to 0pt{\vglue 0mm\line{\small\pfoot\ifnum\count0=\startpage
\copyright\ \gtp\hfill\else
\gt, Volume \thevolumenumber\ (\thevolumeyear)\hfill\fi}\vss
}}
\else
\makeatletter
\def\@oddhead{{\small\lhead\ifnum\count0=\startpage ISSN 1364-0380 (on line)
1465-3060 (printed) \hfill {\lnum\number\count0}\else\ifodd\count0
\def\\{ }\ifx\theshorttitle\relax \thetitle \else\theshorttitle\fi\hfill
{\lnum\number\count0}\else\def\\{ and }{\lnum\number\count0}
\hfill\ifx\theshortauthors\relax 
\theauthors\else\theshortauthors\fi\fi\fi}}\def\@evenhead{\@oddhead}
\def\@oddfoot{\small\lfoot\ifnum\count0=\startpage\copyright\ \gtp\hfill\else
\gt, Volume \thevolumenumber\ (\thevolumeyear)\hfill\fi}
\def\@evenfoot{\@oddfoot}
\makeatother
\fi

\newwrite\gtoutfile
\long\gdef\makeheadfile{  
{\def\\{, }\def\s{ }
\immediate\openout\gtoutfile head.xxx
\immediate\write\gtoutfile{To: math@arxiv.org}
\immediate\write\gtoutfile{Subject: put or rep NNNNN:pppp}
\immediate\write\gtoutfile{--text follows this line--}
\immediate\write\gtoutfile{Proxy-for: \ifx\theasciiauthors\relax
\theauthors\else\theasciiauthors\fi\s<\ifx\theasciiemail\relax\theemail\else\theasciiemail\fi>}
\immediate\write\gtoutfile{\noexpand\\}
\immediate\write\gtoutfile{Authors: \ifx\theasciiauthors\relax
\theauthors\else\theasciiauthors\fi}
{\def\\{ }\immediate\write\gtoutfile{Title: \ifx\theasciititle\relax
\thetitle\else\theasciititle\fi}}
\immediate\write\gtoutfile{Subj-class: GT or SG or MG etc}
\immediate\write\gtoutfile{MSC-class: \theprimaryclass\ifx\thesecondaryclass\relax\else, \thesecondaryclass\fi}
\immediate\write\gtoutfile{Journal-ref: Geom. Topol. \thevolumenumber
(\thevolumeyear) \startpage-\finishpage}
\immediate\write\gtoutfile{Comments: Published by Geometry and Topology at}
\immediate\write\gtoutfile{\s\s http://www.maths.warwick.ac.uk/gt/GTVol\thevolumenumber/paper\thepapernumber.abs.html}
\immediate\write\gtoutfile{\noexpand\\}
\immediate\write\gtoutfile{}
\ifx\theasciiabstract\relax
\immediate\write\gtoutfile{\theabstract}\else
\immediate\write\gtoutfile{\theasciiabstract}\fi
\immediate\write\gtoutfile{}
\immediate\write\gtoutfile{\noexpand\\}
\immediate\write\gtoutfile{}
\immediate\closeout\gtoutfile}}  

\def\maketitlepage{\maketitlep\makeheadfile}
\let\maketitle\maketitlepage

\def\ifplaintex{\expandafter\ifx\csname documentclass\endcsname\relax}

\ifplaintex 
\else
\fi

\expandafter\ifx\csname beginpicture\endcsname\relax
\expandafter\ifx\csname documentclass\endcsname\relax
\input pictex \else
\input prepictex \input pictex \input postpictex \fi\fi

\def\gt{{\mathsurround=0pt\it $\cal G\mskip-2mu$eometry \&\ 
$\cal T\!\!$opology}}        

\def\gtp{{\mathsurround=0pt\it $\cal G\mskip-2mu$eometry \&\ 
$\cal T\!\!$opology $\cal P\!$ublications}}  

\def\lognumber#1{\def\thelognumber{#1}}
\def\volumenumber#1{\def\thevolumenumber{#1}}
\def\papernumber#1{\def\thepapernumber{#1}}
\def\volumeyear#1{\def\thevolumeyear{#1}}

\def\pagenumbers#1#2{\def\startpage{#1}\def\finishpage{#2}}
\def\published#1{\def\publishdate{#1}}
\def\proposed#1{\def\theproposer{#1}}
\def\seconded#1{\def\theseconders{#1}}
\def\received#1{\def\receiveddate{#1}}
\def\revised#1{\def\reviseddate{#1}}
\def\accepted#1{\def\accepteddate{#1}}
\def\asciititle#1{\def\theasciititle{#1}}

\long\def\asciiabstract#1{\long\def\theasciiabstract{#1}}
\def\asciikeywords#1{\def\theasciikeywords{#1}}

\let\\\par\let\thelognumber\relax
\let\thevolumenumber\relax\let\thepapernumber\relax
\let\thevolumeyear\relax\let\thesamplenumber\relax\let\startpage\relax
\let\finishpage\relax\let\publishdate\relax\let\receiveddate\relax
\let\reviseddate\relax\let\accepteddate\relax\let\theasciititle\relax
\let\theasciiauthors\relax
\let\theasciiabstract\relax\let\theasciikeywords\relax
\let\theasciiemail\relax\let\theshortauthors\relax\let\theshorttitle\relax

\long\def\maketitlep{   

\count0=\startpage

\gt\hfill      
\beginpicture
\setcoordinatesystem units <0.33truein, 0.33truein> point at 2.2 0.9
\setplotsymbol ({$\cal G$})
\plotsymbolspacing=9truept
\circulararc 315 degrees from 0 1 center at 0 0
\setplotsymbol ({$\cal T$})
\circulararc 315 degrees from 1 -1 center at 1 0
\endpicture
\break
{\small\ifx\thesamplenumber\relax 
Volume \else Sample
\fi\thevolumenumber\ (\thevolumeyear)
\startpage--\finishpage\nl
Published: \publishdate}
\vglue 0.5truein plus 0.4fil minus 0.1truein

{\parskip=0pt\leftskip 0pt plus 1fil\def\\{\par\smallskip}{\ifplaintex\large
\else\Large\fi\bf\thetitle}\par\medskip}   

\vglue 0pt plus 0.1fil 

{\parskip=0pt\leftskip 0pt plus 1fil\def\\{\par}{\sc\theauthors}
\par\medskip}

\vglue 0pt plus 0.1fil 

{\small\parskip=0pt\let\newline\\
{\leftskip 0pt plus 1fil\def\\{\par}{\sl\theaddress}\par}
\expandafter\ifx\theemail\relax    
\relax\else\vglue 5pt plus 0.02fil minus 2pt\def\\{\stdspace{\rm 
and}\stdspace} 
\cl{Email:\stdspace\tt\theemail}\fi
\ifx\theurl\relax                  
\relax\else\vglue 5pt plus 0.02fil minus 2pt\def\\{\stdspace{\rm 
and}\stdspace}
\cl{URL:\stdspace\tt\theurl}\fi\par}

\vglue 7pt plus 0.3fil minus 3pt

{\bf Abstract}
\vglue 5pt plus 0.1fil minus 2pt

\theabstract

\vglue 7pt plus 0.3fil minus 3pt

{\bf AMS Classification numbers}\quad Primary:\quad \theprimaryclass

Secondary:\quad \thesecondaryclass

\vglue 5pt plus 0.3fil minus 2pt

{\bf Keywords}\quad \thekeywords

\vglue 10pt plus 0.5fil minus 5pt

{\small  Proposed: \theproposer\hfill Received: \receiveddate\nl
Seconded: \theseconders\hfill 
\ifx\reviseddate\relax                         
Accepted: \accepteddate                        
\else
Revised: \reviseddate                          
\fi}
\eject
}       

\let\maketitlepage\maketitlep
\let\maketitle\maketitlepage

\font\phead=cmsl9 scaled 950
\font=cmsl9 scaled 1050
\font\pnum=cmbx10 scaled 913
\font\lnum=cmbx10 
\font\pfoot=cmsl9 scaled 950
\font=cmsl9 scaled 1050
\ifplaintex
\headline{\vbox to 0pt{\vskip -4.5mm\line{\small\phead\ifnum
\count0=\startpage ISSN 1364-0380 (on line)
1465-3060 (printed) \hfill {\pnum\folio}\else\ifodd\count0\def\\{ }%
\ifx\theshorttitle\relax\thetitle\else\theshorttitle\fi\hfill{\pnum\folio}
\else\def\\{ and }{\pnum\folio}\hfill\ifx\theshortauthors\relax\theauthors
\else\theshortauthors\fi\fi\fi}\vss}}
\footline{\vbox to 0pt{\vglue 0mm\line{\small\pfoot\ifnum\count0=\startpage
\copyright\ \gtp\hfill\else
\gt, Volume \thevolumenumber\ (\thevolumeyear)\hfill\fi}\vss
}}
\else
\makeatletter
\def\@oddhead{{\small\lhead\ifnum\count0=\startpage ISSN 1364-0380 (on line)
1465-3060 (printed) \hfill {\lnum\number\count0}\else\ifodd\count0
\def\\{ }\ifx\theshorttitle\relax \thetitle \else\theshorttitle\fi\hfill
{\lnum\number\count0}\else\def\\{ and }{\lnum\number\count0}
\hfill\ifx\theshortauthors\relax 
\theauthors\else\theshortauthors\fi\fi\fi}}\def\@evenhead{\@oddhead}
\def\@oddfoot{\small\lfoot\ifnum\count0=\startpage\copyright\ \gtp\hfill\else
\gt, Volume \thevolumenumber\ (\thevolumeyear)\hfill\fi}
\def\@evenfoot{\@oddfoot}
\makeatother
\fi

\newwrite\gtoutfile
\long\gdef\makeheadfile{  
{\def\\{, }\def\s{ }
\immediate\openout\gtoutfile head.xxx
\immediate\write\gtoutfile{To: math@arxiv.org}
\immediate\write\gtoutfile{Subject: put or rep NNNNN:pppp}
\immediate\write\gtoutfile{--text follows this line--}
\immediate\write\gtoutfile{Proxy-for: \ifx\theasciiauthors\relax
\theauthors\else\theasciiauthors\fi\s<\ifx\theasciiemail\relax\theemail\else\theasciiemail\fi>}
\immediate\write\gtoutfile{\noexpand\\}
\immediate\write\gtoutfile{Authors: \ifx\theasciiauthors\relax
\theauthors\else\theasciiauthors\fi}
{\def\\{ }\immediate\write\gtoutfile{Title: \ifx\theasciititle\relax
\thetitle\else\theasciititle\fi}}
\immediate\write\gtoutfile{Subj-class: GT or SG or MG etc}
\immediate\write\gtoutfile{MSC-class: \theprimaryclass\ifx\thesecondaryclass\relax\else, \thesecondaryclass\fi}
\immediate\write\gtoutfile{Journal-ref: Geom. Topol. \thevolumenumber
(\thevolumeyear) \startpage-\finishpage}
\immediate\write\gtoutfile{Comments: Published by Geometry and Topology at}
\immediate\write\gtoutfile{\s\s http://www.maths.warwick.ac.uk/gt/GTVol\thevolumenumber/paper\thepapernumber.abs.html}
\immediate\write\gtoutfile{\noexpand\\}
\immediate\write\gtoutfile{}
\ifx\theasciiabstract\relax
\immediate\write\gtoutfile{\theabstract}\else
\immediate\write\gtoutfile{\theasciiabstract}\fi
\immediate\write\gtoutfile{}
\immediate\write\gtoutfile{\noexpand\\}
\immediate\write\gtoutfile{}
\immediate\closeout\gtoutfile}}  

\def\maketitlepage{\maketitlep\makeheadfile}
\let\maketitle\maketitlepage

\lognumber{164}
\volumenumber{6}\papernumber{6}\volumeyear{2002}
\pagenumbers{153}{194}
\received{16 February 2001}
\revised{8 July 2001}
\published{30 March 2002}
\accepted{18 March 2002}

\proposed{Cameron Gordon}

\seconded{Joan Birman, Robion Kirby}

\usepackage{amssymb, amsmath} \usepackage{graphicx}
\input labelfig

\theoremstyle{plain} \newtheorem{thm}{Theorem}

\newtheorem{lemma}{Lemma}[section]
\newtheorem{corollary}[lemma]{Corollary}
\newtheorem{proposition}[lemma]{Proposition}

\theoremstyle{definition} \newtheorem{definition}{Definition}[section]

\theoremstyle{remark} \newtheorem{remark}{Remark}[section]
\newtheorem*{notation}{Notation}

\numberwithin{figure}{section}

\nocolon

\begin{document}

\title{Laminar Branched Surfaces in 3--manifolds} 
\asciititle{Laminar Branched Surfaces in 3-manifolds} 

\author{Tao Li}
\address{Department of Mathematics, 401 Math Sciences\\Oklahoma State
University, Stillwater, OK 74078, USA} 

\email{tli@math.okstate.edu}
\url{http://www.math.okstate.edu/\char'176tli}

\begin{abstract}
We define a laminar branched surface to be a branched surface
satisfying the following conditions: (1) Its horizontal boundary is
incompressible; (2) there is no monogon; (3) there is no Reeb
component; (4) there is no sink disk (after eliminating trivial
bubbles in the branched surface).  The first three conditions are
standard in the theory of branched surfaces, and a sink disk is a disk
branch of the branched surface with all branch directions of its
boundary arcs pointing inwards.  We will show in this paper that every
laminar branched surface carries an essential lamination, and any
essential lamination that is not a lamination by planes is carried by
a laminar branched surface.  This implies that a 3--manifold contains
an essential lamination if and only if it contains a laminar branched
surface.
\end{abstract}

\asciiabstract{We define a laminar branched surface to be a branched
surface satisfying the following conditions: (1) Its horizontal
boundary is incompressible; (2) there is no monogon; (3) there is no
Reeb component; (4) there is no sink disk (after eliminating trivial
bubbles in the branched surface).  The first three conditions are
standard in the theory of branched surfaces, and a sink disk is a disk
branch of the branched surface with all branch directions of its
boundary arcs pointing inwards.  We will show in this paper that every
laminar branched surface carries an essential lamination, and any
essential lamination that is not a lamination by planes is carried by
a laminar branched surface.  This implies that a 3-manifold contains
an essential lamination if and only if it contains a laminar branched
surface.}

\asciikeywords{3-manifold, branched surface, lamination}

\maketitlepage

\setcounter{section}{-1}
\section{Introduction}

It has been a long tradition in 3--manifold topology to obtain
topological information using codimension one objects.  Almost all
important topological information has been known for 3--manifolds that
contain incompressible surfaces, eg, \cite{Wa,Th}; other codimension
one objects, such as Reebless foliations and immersed surfaces, have
also been proved fruitful \cite{G1,G2,G3,HRS,N}.  In \cite{GO},
essential laminations were introduced as a generalization of
incompressible surfaces and Reebless foliations and it was proved in
\cite{GO} that if a closed and orientable 3--manifold contains an
essential lamination, then its universal cover is $\mathbb{R}^3$.
More recently, Gabai and Kazez proved that if an orientable and
atoroidal 3--manifold contains a genuine lamination, ie, an essential
lamination that can not be trivially extended to a foliation, then its
fundamental group is negatively curved in the sense of Gromov.

Ever since the invention of essential laminations, branched surfaces
have been a practical tool to study them \cite{GO}.  Gabai and Oertel
have shown that some splitting of any essential lamination is fully
carried by a branched surface satisfying some natural conditions (see
Proposition 1.1) and any lamination carried by such a branched surface
is an essential lamination.  However, these conditions do not
guarantee the existence of essential laminations.  In fact, it was
shown \cite{GO} that even $S^3$ contains branched surfaces satisfying
those conditions.  One of the most important problems in the theory of
essential laminations is to find sufficient conditions for a branched
surface to carry an essential lamination (see Gabai's problem list
\cite{G5}).  In this paper we will show that those standard conditions
in \cite{GO} plus one more, which is that the branched surface does
not contain sink disks, are sufficient and (except for a single
3--manifold) necessary conditions, see section~\ref{S1} for definition
of sink disk.  We call a branched surface satisfying these conditions
a laminar branched surface.

\begin{thm}\label{T01}
Suppose $M$ is a closed and orientable 3--manifold.  Then:
\begin{enumerate}
\item [\rm(a)] Every laminar branched surface in $M$ fully carries an
essential lamination.

\item [\rm(b)] Any essential lamination in $M$ that is not a lamination
by planes is fully carried by a laminar branched surface.
\end{enumerate}
Furthermore, if $\lambda\subset M$ is a lamination by planes (hence
$M=T^3$), then any branched surface carrying $\lambda$ is not a
laminar branched surface.
\end{thm}

Since $T^3=S^1\times S^1\times S^1$ is Haken, and incompressible
surfaces are very special cases of essential laminations, we have:

\begin{thm}\label{T02}
A 3--manifold contains an essential lamination if and only if it
contains a laminar branched surface.
\end{thm}

In many situations, it is easier to construct a branched surface than
to construct an essential lamination.  Theorem~\ref{T01} gives a
criterion to tell whether a branched surface carries an essential
lamination.  It is a very useful theorem.  For example, Delman and Wu
\cite{D,Wu} have shown that many 3--manifolds contain essential
laminations by constructing branched surfaces in certain classes of
knot complements and showing that they carry essential laminations.
Theorem~\ref{T01} can simplify, to some extent, their proofs.  It is
also easy to see that Hatcher's branched surfaces \cite{Ha2} satisfy
our conditions.  Moreover, after splitting the branched surfaces near
the boundary torus such that the train tracks (ie, Hatcher's branched
surfaces restricted to the boundary) become circles, the branched
surfaces also satisfy our conditions (after capping the circles off).
This implies that they are laminar branched surfaces in the manifolds
after the Dehn fillings along these circles.  Hence Hatcher's branched
surfaces carry more laminations than what was shown in \cite{Ha2} and
Theorem~\ref{T01} gives a simpler proof of a theorem of Naimi
\cite{Na}.  More recently, Roberts has constructed taut foliations in
many manifolds using this theorem \cite{R}.

Another interesting question that arose when the concept of essential
lamination was introduced is whether there is a lamination-free theory
for branched surfaces.  In a subsequent paper \cite{L}, we will
discuss this question by proving the following theorem and some
interesting properties of laminar branched surfaces without using
lamination techniques.  Theorem~\ref{T03} is just the branched surface
version of the theorems of Gabai--Oertel \cite{GO} and Gabai--Kazez
\cite{GK}, and it is an immediate corollary of Theorem~\ref{T01}.

\begin{thm}\label{T03}
Let $M$ be a closed and orientable 3--manifold that contains a laminar
branched surface $B$.  Then:
\begin{enumerate}
\item [\rm(i)] The universal cover of $M$ is $\mathbb{R}^3$.

\item [\rm(ii)] If, in addition, the 3--manifold is atoroidal and at least
one component of $M-B$ is not an $I$--bundle, then the fundamental
group of $M$ is word hyperbolic.
\end{enumerate}
\end{thm}

We organize the paper as follows: in section \ref{S1}, we list some
basic definitions and results about essential laminations and give the
definition of laminar branched surfaces; in section 2, we prove some
topological lemmas that we need in the construction of essential
laminations; in sections 3 and 4, we show that every laminar branched
surface carries an essential lamination; in section 5, we prove part
(b) of Theorem~\ref{T01}.

\medskip
\noindent
\textbf{Acknowledgments}\qua  I would like to thank Dave Gabai and Ian Agol
for many very helpful conversations.  I would also like to thank the
referee for many corrections and suggestions.

\section{Preliminaries}\label{S1}

A (codimension one) lamination $\lambda$ in a 3--manifold $M$ is a
foliated, closed subset of $M$, ie, $\lambda$ is covered by a
collection of open sets of the form $\mathbb{R}^2\times\mathbb{R}$
such that, for any open set $U$, $\lambda\cap U=\mathbb{R}^2\times C$,
where $C$ is a closed set in $\mathbb{R}$, and the transition maps
preserve the product structures.  The coordinate neighborhoods of
leaves are of the form $\mathbb{R}^2\times x$ ($x\in C$).

Unless specified, our laminations in this paper are always assumed to
be codimension one laminations in closed and orientable 3--manifolds.
Similar results hold for laminations (with boundary) in 3--manifolds
whose boundary is incompressible.  Let $\lambda$ be a lamination in
$M$, and $M_\lambda$ be the metric completion of the manifold
$M-\lambda$ with the path metric inherited from a Riemannian metric on
$M$.

\begin{definition}\cite{GO}\qua
$\lambda$ is an \textit{essential lamination} in $M$ if it satisfies
the following conditions.
\begin{enumerate}
    \item The inclusion of leaves of the lamination into $M$ induces
    an injection on $\pi_1$.  \item $M_\lambda$ is irreducible.  \item
    $\lambda$ has no sphere leaves.  \item $\lambda$ is
    end-incompressible.
\end{enumerate}
\end{definition}

\begin{figure}[ht!]
\cl{\includegraphics[width=4in]{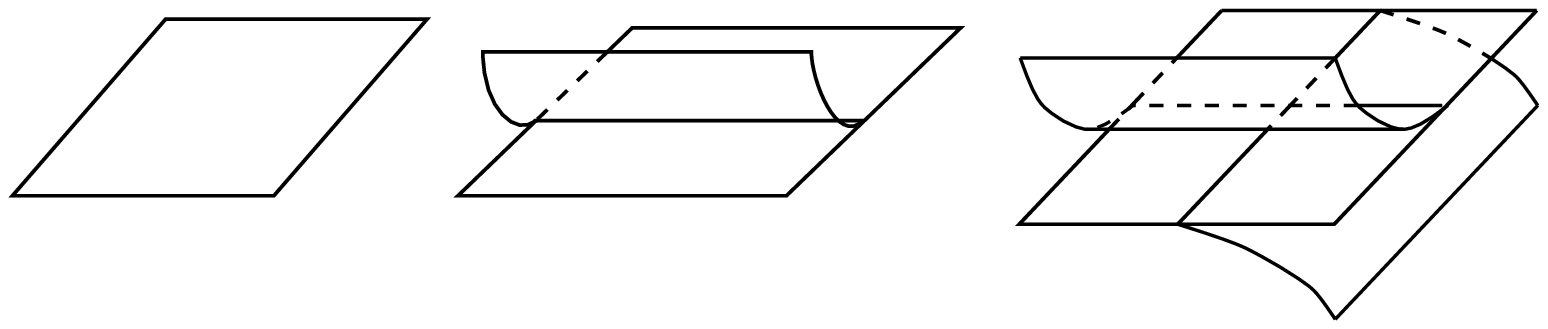}}
\caption{} \label{F11}
\end{figure}

\begin{definition}
A \textit{branched surface} $B$ in $M$ is a union of finitely many
compact smooth surfaces gluing together to form a compact subspace (of
$M$) locally modeled on Figure~\ref{F11}.
\end{definition}

\begin{notation}
Throughout this paper, we denote the interior of $X$ by $int(X)$, and
denote the number of components of $X$ by $|X|$, for any $X$.
\end{notation}

\begin{figure}[ht!]
\vspace{2mm}
\centerline{\small
\SetLabels 
\E(0.2*-0.02){(a)}\\
\E(.7*-0.02){(b)}\\
\E(0.535*1.01){$\partial_vN(B)$}\\
\E(.945*1.04){$\partial_hN(B)$}\\
\endSetLabels 
\AffixLabels{{\includegraphics[width=4in]{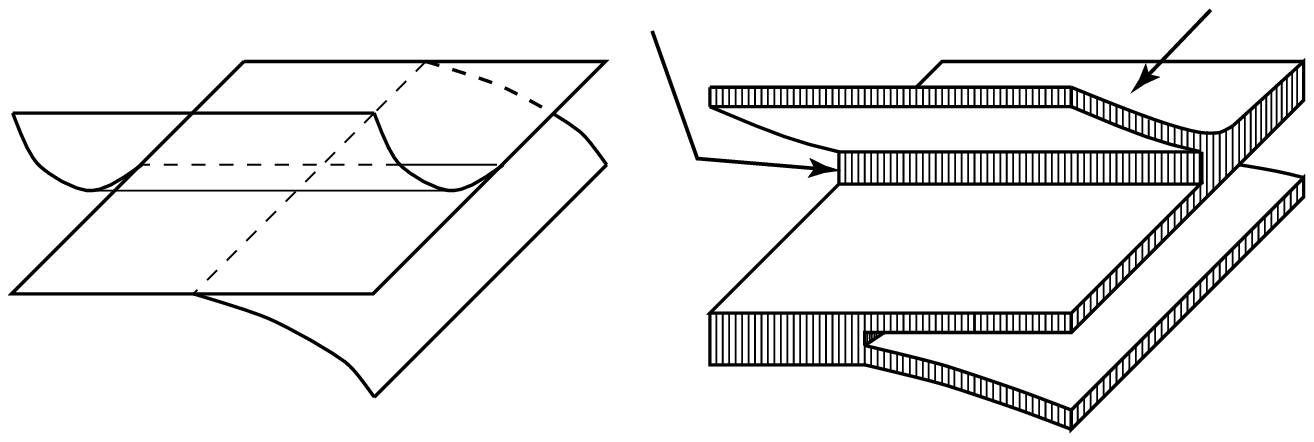}}}}
\vspace{2mm}
\caption{}\label{F12}
\end{figure}

Given a branched surface $B$ embedded in a 3--manifold $M$, we denote
by $N(B)$ a regular neighborhood of $B$, as shown in Figure~\ref{F12}.
One can regard $N(B)$ as an interval bundle over $B$.  We denote by
$\pi \co  N(B)\to B$ the projection that collapses every interval fiber
to a point.  The \textit{branch locus} of $B$ is $L=\{b\in B:$ $b$
does not have a neighborhood homeomorphic to $\mathbb{R}^2 \}$.  So,
$L$ can be considered as a union of smoothly immersed curves in $B$,
and we call a point in $L$ a \textit{double point} of $L$ if any small
neighborhood of this point is modeled on the third picture of
Figure~\ref{F11}.

Let $D_0$ be a component of $B-L$, and $D$ be the closure of $D_0$ in
the path metric (of $B-L$).  Then, $int(D)=D_0$, $\partial D\subset
L$, and those non-smooth points in $\partial D$ are double points of
$L$.  Note that $int(D)=D_0$ is embedded in $B$, but $\partial D$ may
not be embedded in $B$ (there may be two boundary arcs of $D$ that are
glued to the same arc in $L$).  We call $D$ a \textit{branch} of $B$.

The boundary of $N(\!B)$ is a union of two compact surfaces $\partial_h
N(\!B)$ and $\partial_v N(\!B)$.  An interval fiber of $N(B)$ meets
$\partial_h N(B)$ transversely, and intersects $\partial_v N(B)$ (if
at all) in one or two closed intervals in the interior of this fiber.
Note that $\partial_v N(B)$ is a union of annuli, and $\pi (\partial_v
N(B))$ is exactly the branch locus of $B$ (see Figure~\ref{F12}).  We
call $\partial_h N(B)$ the \textit{horizontal boundary} of $N(B)$ and
$\partial_v N(B)$ the \textit{vertical boundary} of $N(B)$.

We say a lamination $\lambda$ is \textit{carried} by $B$ if, after
some splitting, $\lambda$ can be isotoped into $int(N(B))$ so that it
intersects the interval fibers transversely, and we say $\lambda$ is
\textit{fully carried} by $B$ if $\lambda$ intersects every fiber of
$N(B)$.

Gabai and Oertel \cite{GO} found the first relation between essential
laminations and the branched surfaces that carry them.

\begin{figure}[ht!]
\centerline{\small
\SetLabels 
\E(0.2*-0.02){(a)}\\
\E(.8*-0.02){(b)}\\
\endSetLabels 
\AffixLabels{{\includegraphics[width=3.5in]{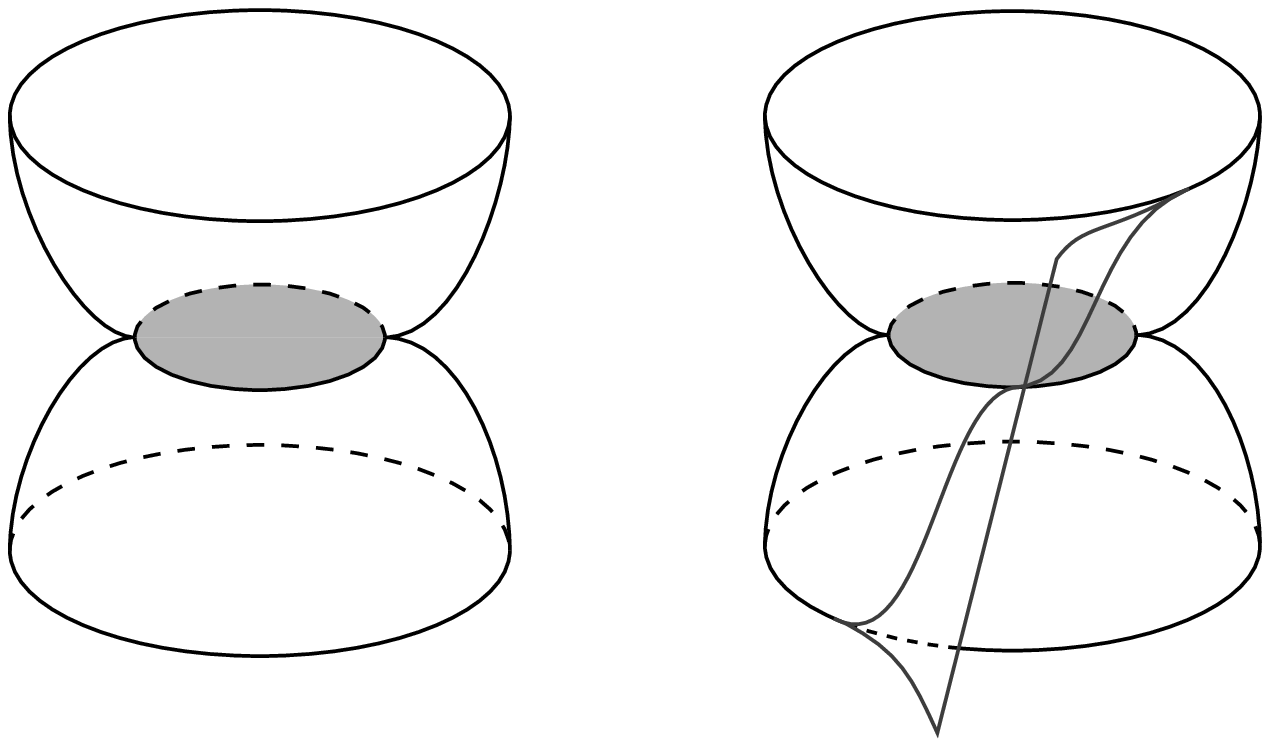}}}}
\vspace{2mm}
\caption{}\label{F13}
\end{figure}

\begin{proposition}[Gabai and Oertel]\label{P11}

{\rm(a)}\qua Every essential lamination is fully carried by a branched surface
with the following properties.
\begin{enumerate}
\item $\partial_h N(B)$ is incompressible in $M-int(N(B))$, no
component of $\partial_h N(B)$ is a sphere, and $M-B$ is irreducible.
\item There is no monogon in $M-int(N(B))$, ie, no disk $D\subset
M-int(N(B))$ with $\partial D=D\cap N(B)=\alpha\cup\beta$, where
$\alpha\subset\partial_v N(B)$ is in an interval fiber of $\partial_v
N(B)$ and $\beta\subset\partial_h N(B)$.
\item There is no Reeb component, ie, $B$ does not carry a torus that
bounds a solid torus in $M$.
\item B has no disk of contact, ie, no disk $D\subset N(B)$ such that
$D$ is transverse to the $I$--fibers of $N(B)$ and $\partial
D\subset\partial_v N(B)$, see Figure~\ref{F13} (a) for an example.

\end{enumerate}
{\rm(b)}\qua If a branched surface with properties above fully carries a
lamination, then it is an essential lamination.
\end{proposition}

However, such branched surfaces may not carry any laminations and they
do not give much information about the 3--manifolds.

\begin{proposition}[Gabai and Oertel]\label{P12}
$S^3$ contains a branched surface satisfying all the conditions in
Proposition~\ref{P11}.
\end{proposition}

It has also been pointed out in \cite{GO} that a twisted disk of
contact is an obvious obstruction for a branched surface to carry a
lamination, because it forces non-trivial holonomy along trivial
curves, which contradicts the Reeb stability theorem (see
Figure~\ref{F13} (b)).

\begin{figure}[ht!]
\centerline{\small
\SetLabels 
\E(0.2*-0.02){(a)}\\
\E(.7*-0.02){(b)}\\
\endSetLabels 
\AffixLabels{{\includegraphics[width=4in]{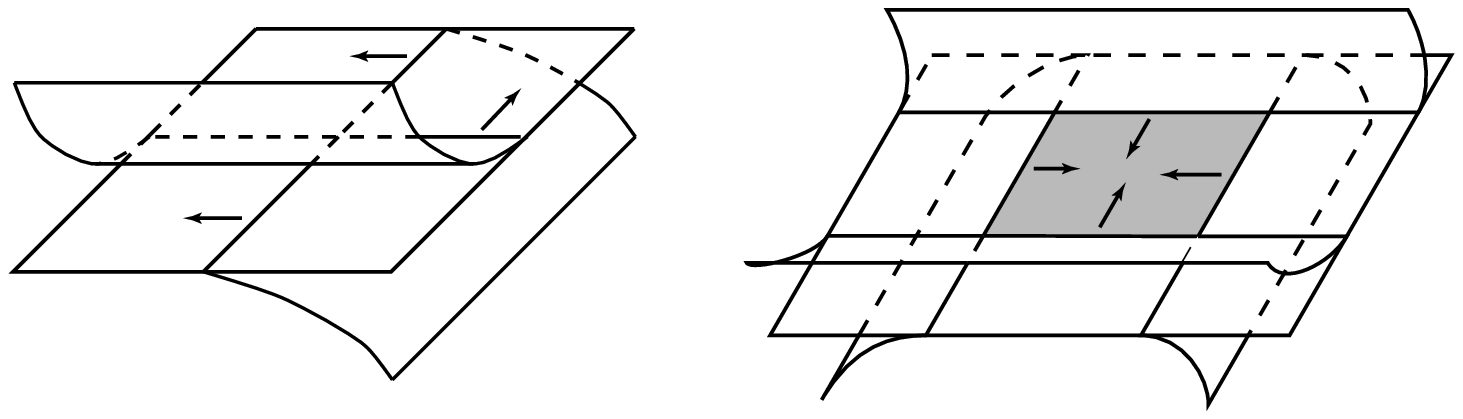}}}}
\vspace{2mm}
\caption{}\label{F15}
\end{figure}

Let $L$ be the branch locus of $B$.  $L$ is a collection of smooth
immersed curves in $B$.  Let $X$ be the union of double points of $L$.
We associate with every component of $L-X$ a vector (in $B$) pointing
in the direction of the cusp, as shown in Figure~\ref{F15}.  We call
it the \textit{branch direction} of this arc.

We call a disk branch of $B$ a \textit{sink disk} if the branch
direction of every smooth arc (or curve) in its boundary points into
the disk.  The standard pictures of disks of contact (Figure~\ref{F13}
(a)) and twisted disks of contact (Figure~\ref{F13} (b)) are all sink
disks by our definition.  Moreover, the disk in Figure~\ref{F15} (b)
is also a sink disk.  Note that a disk of contact can be much more
complicated than Figure~\ref{F13} (a) (see Proposition~\ref{P11} for
the definition of disk of contact).  We will discuss the relation
between a sink disk and a disk of contact in section~\ref{S:facts}
(see Corollary~\ref{C:contact}).

A sink disk can be considered to be a generalized disk of contact.
Here is another way to see this.  In a regular neighborhood of such a
disk, we consider the two components of the complement of $B$ (the one
above the disk and the one below).  The disk is exactly the
intersection of the boundaries of the two components.  Moreover in
each component, one can have a properly embedded disk with smooth
boundary, which is isotopic to the sink disk.

Let $K$ be a component of $M-int(N(B))$.  If $K$ is homeomorphic to a
$3$--ball, then, since $\partial_hN(B)$ is incompressible in
$M-int(N(B))$, $\partial K$ consists of two disk components of
$\partial_hN(B)$ and an annulus component of $\partial_vN(B)$.
Moreover, we can give $K$ a fiber structure $D^2\times I$, with
$D^2\times\partial I\subset\partial_hN(B)$ and $\partial D^2\times
I\subset\partial_vN(B)$.  We call $K$ a ${D^2\times I}$
\textit{region} in $M-int(N(B))$, $D^2\times\partial I$ the
\textit{horizontal boundary} of $K$ and $\partial D^2\times I$ the
\textit{vertical boundary} of $K$.

\begin{definition}
Let $D_1$ and $D_2$ be the two disk components of the horizontal
boundary of a $D^2\times I$ region $K$ in $M-int(N(B))$.  Hence, $D_1$
and $D_2$ are also two disk components of $\partial_hN(B)$ and
$D_1\cup D_2=D^2\times\partial I$.  Thus, $\pi(\partial
D_1)=\pi(\partial D_2)$ is a circle in the branch locus $L$, where
$\pi\co N(B)\to B$ is the collapsing map.  If $\pi$ restricted to the
interior of $D_1\cup D_2$ is injective, ie, the intersection of any
$I$--fiber of $N(B)$ with $int(D_1)\cup int(D_2)$ is either empty or a
single point, then we call $K$ a \textit{trivial} ${D^2\times
I}$ \textit{region}, and we say that $\pi(D_1\cup D_2)$ forms a
\textit{trivial bubble} in $B$.
\end{definition}

Let $K=D^2\times I$ be a trivial $D^2\times I$ region.  Then, after
collapsing each $I$--fiber of $K=D^2\times I$ to a point, $N(B)\cup K$
becomes a fibered neighborhood of another branched surface with the
induced fiber structure from $N(B)$.  Thus, if $B$ contains a trivial
bubble, we can pinch $B$ to get another branched surface by collapsing
the $I$--fibers in the corresponding trivial $D^2\times I$ region, and
the new branched surface after this pinching preserves the properties
1--4 in Proposition \ref{P11}.  There is really no difference between
the branched surface before this pinching and the one after the
pinching.  It is easy to see that a branched surface carries a
lamination if and only if, after we collapse all trivial bubbles in
$B$ as above, the new branched surface carries a lamination.  Not all
$D^2\times I$ regions are trivial, eg, we cannot collapse all the
$D^2\times I$ regions in a standard Reeb component.  Moreover, if we
blow an ``air bubble" into the interior of a sink disk, it will
destroy the sink disk by definition but nothing really changes.  So,
in this paper, we always assume $B$ contains no trivial bubble.

\begin{definition}
A branched surface $B$ in $M$ is called a \textit{laminar branched
surface} if it satisfies conditions 1--3 in Proposition~\ref{P11}, and
$B$ has no sink disk (after we collapse all the trivial bubbles as
described above).
\end{definition}

In this paper, we will also use some techniques about train tracks.
We refer readers to \cite{PH} section 1.1 for basic definitions and
properties about train tracks.  Let $D$ be a disk and $\tau$ be a
train track in $D$.  Suppose $W$ is a closed disk embedded in the
plane whose boundary is piecewise smooth with $k\ge 0$ discontinuities
in the tangent.  Let $h\co W\to D$ be a $C^1$ immersion which is an
embedding of the interior of $W$, and $h(W)\subset\tau$.  Let
$\Delta=h(W)$, and we denote the image (under $h$) of $int(W)$ by
$int(\Delta)$. Note that $int(\Delta)$ is an embedded open disk in
$D$, but $\Delta-int(\Delta)$ may not be embedded in $D$.  We call
$\Delta$ a $k$--\textit{gon}, if $k>0$ and $\tau-int(\Delta)$ is a sub
train track of $\tau$ in $D$.  So, $\Delta-int(\Delta)$ is an immersed
circle with $k$ prongs, each smooth arc in $\Delta-int(\Delta)$ is
carried by $\tau$, and the $k$ non-smooth points in $\partial\Delta$
are switches (non-manifold points) of $\tau$.  If $k=1$, we also call
$\Delta$ a monogon, and if $k=2$, we also call $\Delta$ a bigon.  If
$\tau\cap int(\Delta)=\emptyset$, we call $\Delta$ a $k$--\textit{gon
component} of $D-\tau$.  We call $\Delta$ a \textit{smooth disk} if
$k=0$ and $\tau-int(\Delta)$ is a sub train track of $\tau$ in $D$,
ie, $int(\Delta)$ is an embedded disk, and $\Delta-int(\Delta)$ is a
circle carried by $\tau$.  We call $\Delta$ a \textit{smooth disk
component} of $D-\tau$ if $\tau\cap int(\Delta)=\emptyset$.  Let
$N(\tau)$ be a fibered neighborhood of $\tau$.  Then, the $\Delta$
above corresponds to an embedded disk in $D$; if $\Delta$ is a smooth
disk, $\partial\Delta$ corresponds to an embedded circle in $N(\tau)$
transversely intersecting the $I$--fibers of $N(\tau)$; if $\Delta$ is
a $k$--gon, $\partial\Delta$ corresponds to an embedded circle in
$N(\tau)$ consisting of $k$ arcs, each of which is transverse to the
$I$--fibers of $N(\tau)$.  Throughout this paper, when we talk about an
object in $D$ with respect to the train track $\tau$, we
simultaneously use the same notation to denote the corresponding
object in $D$ with respect to $N(\tau)$.

Let $\Delta$ be a $k$--gon as above.  We call $\Delta-int(\Delta)$
(ie, $h(W-int(W))$ the boundary of $\Delta$, which we denote by
$\partial\Delta$.  We call the image (under $h$) of a non-smooth point
of $\partial W$ a \textit{vertex} of the $k$--gon $\Delta$, and call
the image (under $h$) of a smooth arc between two non-smooth points in
$\partial W$ an \textit{edge} of the $k$--gon $\Delta$.

\section{Some topological lemmas}\label{S:facts}
In this section, we explore topological and combinatorial properties
of laminar branched surfaces by proving some lemmas.  Lemmas
\ref{L:dl} and \ref{L:ba} will be used in section~\ref{S:construction}
to guarantee that a part of the lamination constructed in
section~\ref{S:construction} satisfies a technical condition in a
lemma.  Lemma~\ref{L:B-} is interesting in its own right.  In
particular, we prove Corollary~\ref{C:contact}, which basically says
that the condition of no sink disks implies that there is no disk of
contact.  Note that the condition of no disks of contact plays an
important role in the proof of Proposition~\ref{P11} (b) \cite{GO}.

\begin{figure}[ht!]
\vspace{2mm}
\centerline{\small
\SetLabels 
\E(0.2*-0.02){$B$}\\
\E(.7*-0.02){$B-S$}\\
\E(0.25*.4){$S$}\\
\endSetLabels 
\AffixLabels{{\includegraphics[width=4in]{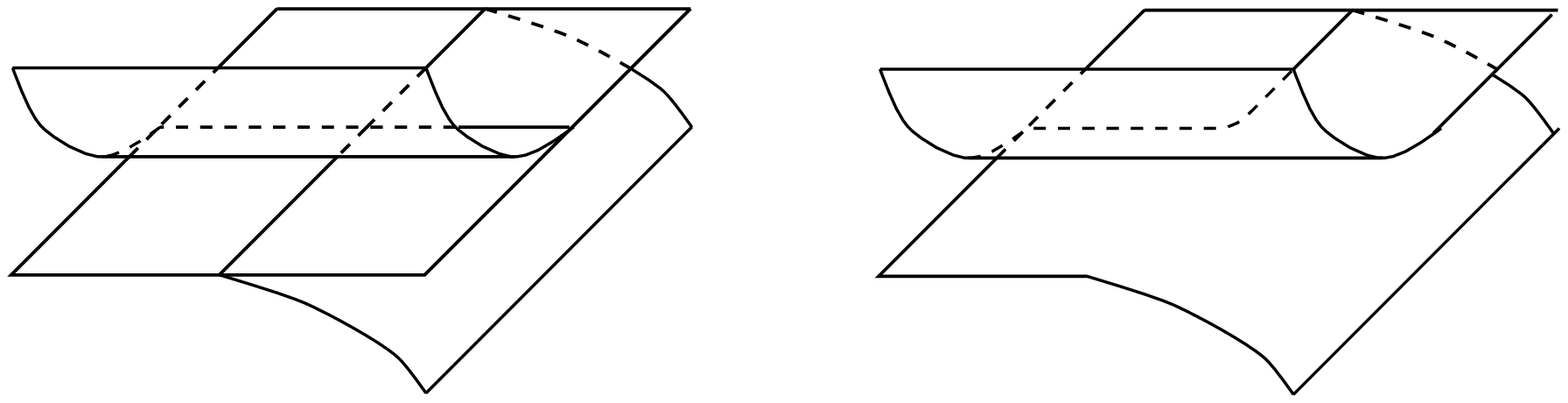}}}}
\vspace{2mm}
\caption{}\label{delete}
\end{figure}

Let $B$ be a laminar branched surface, and $S$ be a branch of $B$.
The boundary of $S$ is piecewise smooth, and each smooth arc in
$\partial S$ has a transverse direction induced from the branch
direction of the corresponding arc in $L$, where $L$ denotes the
branch locus throughout this paper.  Then, we can consider $B$ to be
the object obtained by gluing all the branches of $B$ together along
their boundaries according to the branch directions.  If the branch
direction of each smooth arc in $\partial S$ points out of $S$ and
there are no two arcs in $\partial S$ glued together (to the same arc
in $L$), then $B-int(S)$ naturally forms another branched surface, as
shown in Figure \ref{delete}.  We denote this branched surface
($B-int(S)$) by $B^-$.  Note that if two arcs in $\partial S$ are
identified to the same arc in $L$, $B-int(S)$ is not a branched
surface anymore near this arc.  Moreover, no three arcs in $\partial
S$ can be identified to the same arc in $L$, because otherwise, one of
the three arcs must have induced direction (from the branch direction)
pointing into $S$.

\begin{definition}
Let $S$ be a disk branch of $B$ with branch direction of each boundary
edge pointing out of $S$.  If there are no two arcs in the boundary of
$S$ identified to the same arc in $L$, we call $S$ a \textit{removable
disk}.  If $B$ contains no removable disk, we say that $B$ is
\textit{efficient}.
\end{definition}

\begin{lemma}\label{L:B-}
Let $B$ be a laminar branched surface and $S$ be a removable disk in
$B$.  Then, $B^-=B-int(S)$ is also a laminar branched surface.
\end{lemma}
\begin{proof}
We first note that we have assumed our laminar branched surface $B$
does not have any trivial bubble.  Then, $B^-$ does not contain any
trivial bubble either, since $B$ can be considered as the branched
surface obtained by adding a branch $S$ to $B^-$ and if we add a disk
branch inside a trivial bubble of $B^-$, we always get a trivial
bubble in $B$.

Now, we show that $B^-$ has no sink disk.  Suppose $D$ is a sink disk
in $B^-$, ie, $D$ is a disk branch with branch directions of its
boundary arcs pointing into $D$.  If $D\cap\partial S$ is a union of
arcs, then $D$ is cut into pieces by $\partial S$, but at least one of
these pieces is a sink disk in $B$, which gives a contradiction.  If
$D\cap\partial S$ contains a circle, since $S$ is a disk branch of $B$
and $\partial_hN(B)$ has no sphere component, $\partial S$ must be a
circle that bounds a disk $D'$ in $D$ and the branch direction of
$\partial S$ must point out of $D'$.  Thus, $S\cup D'$ forms a trivial
bubble in $B$, as $M-B$ is irreducible, which contradicts our
assumption of no trivial bubbles.

Since any surface (or lamination) carried by $B^-$ must also be
carried by $B$, $B^-$ has no Reeb component.  Since no component of
$\partial_hN(B)$ is a 2--sphere, it is easy to see that no component of
$\partial_hN(B^-)$ is a 2--sphere.  Moreover, $M-B^-$ is irreducible,
since a reducing sphere intersects $S$ in loops, which bound disks in
the disk branch $S$, and the irreducibility follows from a standard
cut and paste argument.  So, we only need to show that
$\partial_hN(B^-)$ is incompressible in $M-int(N(B^-))$, and there is
no monogon in $M-B^-$.

Note that $\partial_vN(B^-)$ has a natural fiber structure with every
$I$--fiber a subarc of an $I$--fiber of $N(B^-)$.  Let $\psi \co N(B^-)\to
N_{B^-}$ be a map such that:

\begin{enumerate}
\item $\psi$ collapses every $I$--fiber of $\partial_vN(B^-)$ to a
point;
\item $\psi$, when restricted to $int(N(B))$ and
$int(\partial_hN(B))$, is a homeomorphism.
\end{enumerate}  
Figure \ref{collapse} is a schematic picture of $\psi$.  We denote the
image of $\psi$ by $N_{B^-}$.  Let $\partial_hN_{B^-}$ be the image of
$int(\partial_hN(B^-))$ (under the map $\psi$).  If $\partial_hN(B^-)$
is compressible in $M-int(N(B^-))$, then $\partial_hN_{B^-}$ is
compressible in $M-int(N_{B^-})$.

\begin{figure}[ht!]
\vspace{2mm}
\centerline{\small
\SetLabels 
\E(0.2*-0.15){$N(B^-)$}\\
\E(.8*-0.15){$N_{B^-}$}\\
\E(0.5*.6){$\psi$}\\
\endSetLabels 
\AffixLabels{{\includegraphics[width=4in]{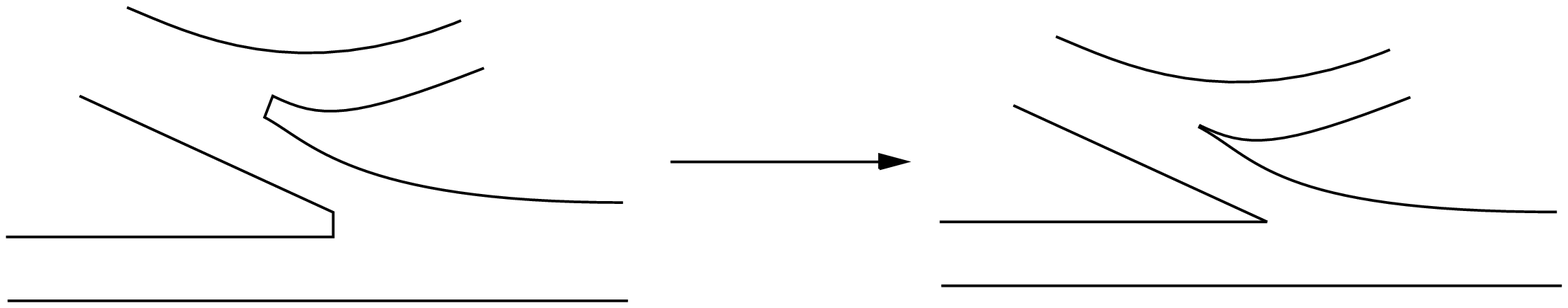}}}}
\vspace{2mm}
\caption{}\label{collapse}
\end{figure}

The component $S$ is a surface with $\partial S\subset B^-$ and
$int(S)$ embedded in $M-B^-$.  The branched surface $B$ can be
considered as the union of $B^-$ and $S$ by smoothing out $\partial S$
according to the branch direction.  In the same way as adding $S$ to
$B^-$, we can add $S$ to $N_{B^-}$.  We can view $S$ as a surface
properly embedded in $M-int(N_{B^-})$ with $\partial S$ piecewise
smoothed out according to the branch direction.  We consider this
complex $N_{B^-}\cup S$.  Note that if we collapse every $I$--fiber of
$N_{B^-}$ to a point, $N_{B^-}\cup S$ becomes $B$; and if we thicken
$N_{B^-}\cup S$ a little, it becomes $N(B)$.  Let $E$ be a compressing
disk in $M-int(N_{B^-})$ with $\partial E\subset\partial_hN_{B^-}$.
We may assume that the compressing disk $E$ intersects $S$
transversely except at $\partial E\cap\partial S$.  We also assume
$|E\cap S|$ (the number of components of $E\cap S$) is minimal among
all compressing disks for $\partial_hN_{B^-}$.  Thus, $E\cap S$
contains no closed circles, otherwise, since a circle of intersection
bounds a disk in $S$, a standard cutting and pasting argument gives us
a compressing disk with fewer intersection curves (with $S$).

Next, we show that $E\cap S\ne\emptyset$.  Suppose $E\cap
S=\emptyset$.  Let $K$ be the closure of the component of $M-N_{B^-}$
that contains $S$.  Then, $E\subset K$, otherwise, it contradicts the
assumption that $\partial_hN(B)$ is incompressible in $M-int(N(B))$.
As $E\cap S=\emptyset$, we can simultaneously consider $E$ to be a
disk embedded in $M-int(N(B))$ with $\partial E$ a smooth nontrivial
circle in $\partial_hN(B)$.  Since $\partial_hN(B)$ is incompressible,
there is an embedded disk $E'\subset\partial_hN(B)$ with $\partial
E=\partial E'$ and $E\cup E'$ bounds an embedded 3--ball in
$M-int(N(B))$.  There are a pair of disks in $\partial_hN(B)$, which
we denoted by $S_1$ and $S_2$, such that $\pi(S_1)=\pi(S_2)=S$ ($S_1$
and $S_2$ can be considered as two sides of $S$).  For each smooth arc
$\alpha\subset\partial S$, let $\alpha_1\subset\partial S_1$ and
$\alpha_2\subset\partial S_2$ be the two corresponding arcs such that
$\pi(\alpha_1)=\pi(\alpha_2)=\alpha$.  Then, since the branch
direction of $\partial S$ points out of $S$, either $\alpha_1$ or
$\alpha_2$ must lie in the boundary of $\partial_hN(B)$
(ie, $\partial_hN(B)\cap\partial_vN(B)$).  Thus, $E'\subset
int(\partial_hN(B))$ cannot intersect both $S_1$ and $S_2$.
Therefore, there must be a smooth disk $\Delta$ embedded in
$\partial_hN_{B^-}\cup S$ such that $\partial\Delta=\partial E$ and
$\Delta\cup E$ bounds an embedded 3--ball in $K$.  Since
$\partial_hN(B)$ is incompressible, $\Delta\cap S\ne\emptyset$.  Note
that the purpose of the argument about $S_1$ and $S_2$ is to show that
$\Delta$ cannot cover $S$ from both sides of $S$ and hence the 3--ball
bounded by $E\cup\Delta$ is embedded.  Since $E\cap S=\emptyset$, $S$
must lie in the interior of $\Delta$.  Since $S$ is a disk and since
$\partial_hN(B)$ is incompressible in $M-int(N(B))$, $K$ must be a
solid torus with $\partial K\subset\partial_hN_{B^-}$, and $K\cup S$
forms a Reeb component, which contradicts our assumptions on $B$.

So, $E\cap S\ne\emptyset$.  Since $E$ and $S$ are properly embedded in
$M-int(N_{B^-})$, $E\cap S$ is a union of disjoint simple arcs in $E$.
Since $\partial S$ is smoothed out according to the branch direction,
the union of $\partial E$ and $E\cap S$ is a train track in $E$ with
all the switches (ie non-manifold points) in $\partial E$.
Moreover, since $M-B$ has no monogon, $E-E\cap S$ has no monogon
component.  Hence, by a standard index argument, $E-E\cap S$ must have
a smooth disk component (see section \ref{S1} for our definitions of
smooth disk component and smooth disk).  We denote this smooth disk
component of $E-E\cap S$ by $E_0$.  Since $E\cap S$ is a union of
disjoint properly embedded arcs, $E-E_0$ is a union of bigons, which
we denote by $E_1, \dots,E_n$ ($E_i$ may contain other components of
$E\cap S$).  The boundary of each bigon $E_i$ consists of two edges,
$\alpha_i$ and $\beta_i$, where $\alpha_i=E_i\cap E_0$ and
$\beta_i=E_i\cap\partial E$.

Since $\partial_hN(B)$ is incompressible, $E_0$ must be parallel to
$\partial_hN_{B^-}\cup S$, ie, there is a smooth disk $\Delta$ in
$\partial_hN_{B^-}\cup S$ such that $\partial\Delta=\partial E_0$ and
$\Delta\cup E_0$ bounds a 3--ball.  Note that by the argument about
$S_1$ and $S_2$ above, only one side of $S$ (near $\partial S$) can be
in the interior of a smooth surface that corresponds to
$\partial_hN(B)$.  Hence, $\Delta$ must be embedded in
$\partial_hN_{B^-}\cup S$, and $\Delta\cup E_0$ bounds an embedded
3--ball $T$.  If $(E-E_0)\cap T\ne\emptyset$, then (since $E$ is
embedded) there must be another smooth disk component $E_0'$ (of
$E-E\cap S$) lying in $(E-E_0)\cap T$, and there is a sub-disk of
$\Delta$, say $\Delta'$, such that $\partial E_0'=\partial\Delta'$ and
$E_0'\cup\Delta'$ bounds a $3$--ball inside $T$.  Thus, by choosing an
appropriate smooth disk component, we can assume that $(E-E_0)\cap
T=\emptyset$.

Since there are no two arcs in $\partial S$ identified to the same arc
in $L$, $\partial S$ is embedded in $\partial N_{B^-}$.  Hence,
$\Delta\cap\partial S$ is a union of disjoint curves that are properly
embedded in $\Delta$.  Since $S$ is a disk, $\Delta\cap\partial S$
contains no closed curves, and hence $\Delta\cap S$ is a union of
disks in $\Delta$.  We denote the components of $\Delta\cap S$ by
$F_1,\dots,F_m$.  Then, each arc in $\partial F_i-\partial S$ is one
of the $\alpha_j$'s (in $\partial E_j$'s) defined before.

Let $\hat{F}_i$ be the union of $F_i$ and those $E_j$'s that share
boundary edges $\alpha_j$'s with $F_i$.  By our construction,
$\cup_{i=1}^n\alpha_i\subset\cup_{i=1}^mF_i$, and hence
$\cup_{i=1}^nE_i\subset\cup_{i=1}^m\hat{F}_i$.  Since each $E_j$ is a
bigon with $\partial E_j=\alpha_j\cup\beta_j$ ($\beta_j\subset\partial
E-\partial E_0\subset\partial_hN_{B^-}$), and since $(E-E_0)\cap
T=\emptyset$, $\hat{F}_i$ is an embedded disk with
$\partial\hat{F}_i\subset\partial_hN_{B^-}$.  Moreover, after pushing
$\hat{F}_i$ out of $T$, $\hat{F}_i\cap S$ has fewer components than
$E\cap S$.  Since we have assumed that $E\cap S$ has the least number
of components among compressing disks, $\hat{F}_i$ cannot be a
compressing disk.  So, $\hat{F}_i$ can be homotoped into
$\partial_hN_{B^-}$ fixing $\partial\hat{F}_i$, for any $i$.  However,
we can then first homotope $E_0$ into $\Delta$ fixing $\partial E_0$
and each $E_i$, then we homotope every $\hat{F}_i$ into
$\partial_hN_{B^-}$ fixing $\partial\hat{F}_i$.  Since $\Delta
-\partial_hN_{B^-}=\cup_{i=1}^mF_i$, and since
$\cup_{i=1}^nE_i\subset\cup_{i=1}^m\hat{F}_i$, after those homotopies
above, we have homotoped $E$ into $\partial_hN_{B^-}$ fixing $\partial
E$, which contradicts the assumption that $E$ is a compressing disk.
Therefore, $\partial_hN_{B^-}$ is incompressible in $M-int(N_{B^-})$,
and hence $\partial_hN(B^-)$ is incompressible in $M-int(N(B^-))$.

Using a similar argument, we can show that $M-N_{B^-}$ contains no
monogons.  Note that since the argument in this case is very similar
to the one above, we keep the same notation, and refer many details to
the argument above.  Now, we let $E$ be a monogon, and suppose $E$
intersects $S$ transversely except at $\partial S$.  We assume $E\cap
S$ has the least number of components among all monogons.  Note that
$E\cap S\ne\emptyset$, since $M-B$ contains no monogons.  As in the
argument above, $(E\cap S)\cup\partial E$ is a train track in $E$.
Since $M-B$ contains no monogons, $E-E\cap S$ contains no monogon
component.  Hence, by a standard index argument, there must be a
smooth disk component in $E-E\cap S$, which we denote by $E_0$.  In
this case, $E-E_0$ is a union of bigons and one $3$--gon (ie, a disk
with three prongs).  Let $E_1, \dots, E_n$ be the components of
$E-E_0$, and suppose $E_1$ is the disk with three prongs.  As before,
there is a smooth disk $\Delta$ in $\partial_h N_{B^-}\cup S$ such
that $\partial\Delta=\partial E_0$ and $\Delta\cup E_0$ bounds an
embedded $3$--ball.  We can define $F_i$'s and $\hat{F}_i$'s as above.
However, in this case, the $\hat{F}_k$ that contains $E_1$ must be a
monogon.  After pushing this $\hat{F}_k$ out of the $3$--ball bounded
by $\Delta\cup E_0$, we get a monogon with fewer intersection curves
(with $S$), which gives a contradiction.
\end{proof}

\begin{remark}\label{R:construction}
\begin{enumerate}
\item Let $B$ and $B^-$ be as above.  By Corollary~\ref{C:contact},
$B$ and $B^-$ have no disks of contact.  Hence, if $B^-$ fully carries
a lamination, using the techniques in \cite{G4,G1} (see also section
\ref{S:extend}), we can construct a lamination fully carried by $B$.
Then, by Proposition \ref{P11} (b), this lamination is an essential
lamination.

\item Let $B_1, \dots, B_m$ be a series of branched surfaces, $L_i$ be
the branch locus of $B_i$ (for any $i$), and $S_i$ ($i<m$) be a
removable disk of $B_i$.  Suppose $B_{i+1}=B_i-int(S_i)$ ($i<m$).  If
$B_1$ is a laminar branched surface, then by Lemma \ref{L:B-}, we can
inductively show that each $B_i$ is a laminar branched surface.
Moreover, as we point out above, if $B_m$ fully carries a lamination,
we can inductively construct a lamination for each $B_i$.  For any
laminar branched surface $B=B_1$, there always exist such a series of
branched surfaces such that $B_m$ is efficient.  If we can construct a
lamination carried by $B_m$, we can inductively extend this lamination
to a lamination fully carried by $B=B_1$.

\item Although we used the hypothesis that $S$ is a disk in the proof,
Lemma~\ref{L:B-} is still true if $S$ is not a disk.
\end{enumerate} 
\end{remark}

\begin{definition}\label{D:parallel}
Let $S_1$ and $S_2$ be two surfaces or arcs in $N(B)$ that are
transverse to the $I$--fibers of $N(B)$.  We say that $S_1$ and $S_2$
are \textit{parallel} if there is an embedding $H\co  S\times [1,2]\to
N(B)$ such that $H(S\times\{i\})=S_i$ ($i=1,2$) and
$H(\{x\}\times[1,2])$ is a subarc of an $I$--fiber of $N(B)$ for any
$x\in S$.
\end{definition}

\begin{definition}
Let $B$ be a branched surface and $D$ be an embedded disk in $N(B)$
that is transverse to the $I$--fibers of $N(B)$.  Suppose $\partial
D\subset\pi^{-1}(L)$, where $L$ is the branch locus of $B$.  Then,
every arc in $\partial D$ has an induced direction that is consistent
with the branch direction of the corresponding arc in $L$.  We call
$D$ a \textit{generalized sink disk} if the induced direction of every
arc in $\partial D$ points into $D$.  Note that if $\pi^{-1}(L)\cap
int(D)=\emptyset$, $\pi(D)$ is a sink disk.
\end{definition}

\begin{lemma}\label{L:gsd}
Let $B$ be a laminar branched surface.  Then, $N(B)$ contains no
generalized sink disk. \end{lemma}
\begin{proof}
Suppose $D$ is a generalized sink disk.  We first show that there must
be a subdisk of $D$, which we denote by $D'$, such that $D'$ is a
generalized sink disk, and $\pi|_{D'}$ is injective (ie, the
intersection of each $I$--fiber of $N(B)$ with $D'$ is either empty or
a single point).

Let $n$ be the maximal number of intersection points of $D$ with any
$I$--fiber of $N(B)$, and $X_n$ be the union of $I$--fibers of $N(B)$
whose intersection with $D$ consists of $n$ points.  We assume $n>1$,
otherwise, $D'=D$.  We use induction on $n$.  Since the induced
direction of every arc in $\partial D$ points into $D$, $X_n\cap D$ is
a collection of compact subsurfaces of $D$.  Moreover, since $n$ is
maximal, the boundary of $X_n\cap D$ lies in $\pi^{-1}(L)$ with
direction (induced from the branch direction of $L$) pointing into
$X_n\cap D$.  Let $P_1,\dots,P_k$ be the components of $X_n\cap D$,
and hence each $P_i$ is a planar surface in $D$.  We call a boundary
circle of $P_i$ the \emph{outer boundary} of $P_i$ if it bounds a disk
in $D$ that contains $P_i$.  We denote the outer boundary of $P_i$ by
$\alpha_i$ ($i=1,\dots k$) and let $D_i$ be the disk bounded by
$\alpha_i$ in $D$.  Hence, $P_i\subset D_i$.  Without loss of
generality, we can assume that $\alpha_1$ is an inner most circle
(among the $\alpha_i$'s), ie, $D_i\not\subset D_1$ for any $i\ne 1$.
Then, $\partial D_1\subset\pi^{-1}(L)$ and the induced direction of
every arc in $\partial D_1$ points into $D_1$.  Hence, $D_1$ is a
generalized sink disk.  Next, we show that we can assume the maximal
number of intersection points of $P_1$ with any $I$--fiber of $N(B)$ is
less than $n$.

Note that $X_n$ can be considered as an $I$--bundle over a compact
surface (if one collapses every $I$--fiber in $X_n$ to a point, since $n$
is maximal, one does not get any branching).  Thus, if $P_1$ contains
all $n$ points of the intersection, $X_n$ must be a twisted $I$--bundle
over a nonorientable surface, $n=2$, and $P_1$ double covers a
nonorientable surface.  If two inner boundary components, say $c_1$
and $c_2$, of $P_1$ are two parallel curves in a vertical boundary
component of the $I$--bundle $X_n$, we can replace the disk (in $D$)
bounded by $c_1$ by a disk that is parallel to the disk (in $D$)
bounded by $c_2$.  By this cutting and pasting, we can assume that
$P_1$ is an annulus that double covers a M\"{o}bius band.  Moreover,
the outer boundary $\alpha_1$ of $P_1$ bounds a disk that is parallel
to the disk (in $D$) bounded by the inner boundary of $P_1$.  By
capping $\alpha_1$ off using this disk, we get an embedded 2--sphere
that double covers a projective plane carried by $B$.  Then, by
applying the train track argument below to this 2--sphere, one either
gets a generalized sink disk $D'$ in this 2--sphere with $\pi|_{D'}$
injective, or gets a removable disk disjoint from a generalized sink
disk and eventually has a contradiction similar to the argument below.
Note that if we assume $M$ to be irreducible here, $M$ must be
$\mathbb{R}P^3$ and it is easy to conclude that this case cannot
happen.  Hence, we may assume the maximal number of intersection
points of $P_1$ with any $I$--fiber of $N(B)$ is less than $n$.

Since $\alpha_1$ is innermost, the maximal number of intersection
points of $D_1$ with any $I$--fiber of $N(B)$ must be less than $n$.
Inductively, we can eventually find a generalized sink disk $D'\subset
D$ such that $\pi|_{D'}$ is injective.  Note that $\pi(D')$ is not
necessarily a sink disk by definition because $\pi^{-1}(L)\cap
int(D')$ may not be empty.

In the remaining part of the proof, we will show that there is a
removable disk $S$ in $B$ such that the fibered neighborhood of the
branched surface $B^-=B-int(S)$ also contains a generalized sink disk.
Let $D'$ be a generalized sink disk such that $\pi|_{D'}$ is
injective.  Moreover, we may assume that $D'$ contains no subdisk that
is a generalized sink disk.  Note that $\pi^{-1}(L)\cap
int(D')\ne\emptyset$, otherwise, $\pi(D')$ is a sink disk which
contradicts the hypothesis that $B$ is a laminar branched surface.

We fix a normal direction for $D'$.  For every point $x\in int(D')$,
let $J_x$ be the $I$--fiber of $N(B)$ that contains $x$.  Then, $J_x-x$
has two components.  According to the fixed normal direction of $D'$,
we say that the points in one component of $J_x-x$ are on the positive
side of $x$, and points in the other component of $J_x-x$ are on the
negative side of $x$.  Let $G$ be the union of $x\in D'$ such that
$J_x\subset\pi^{-1}(L)$ and $\partial_vN(B)\cap J_x$ contains a
component on the positive side of $x$.  Note that if $\pi(J_x)$ is a
double point of $L$, $\partial_vN(B)\cap J_x$ consists of two disjoint
arcs.  Then, by the construction of $G$ and the local model
(Figure~\ref{F12}) of a branched surface, $G\cup\partial D'$ is a
trivalent graph and each edge has a direction induced from the branch
direction.  As shown in Figure \ref{traintrack}, this trivalent graph
$G\cup\partial D'$ can be deformed into a transversely oriented train
track $\tau$ according to the directions of the edges in $G$.  Since
the direction of every arc in $\partial D'$ points into $D'$ and
$\tau$ is transversely oriented, $\partial D'$ is a smooth circle in
$\tau$.  Note that, by choosing an appropriate normal direction for
$D'$, we can assume $G\ne\emptyset$, since $\pi^{-1}(L)\cap
int(D')\ne\emptyset$.  By a standard index argument, $D'-\tau$ must
have a smooth disk component, ie, there is a smooth circle in $\tau$
which is the boundary of the closure of a disk component of $D'-\tau$.
We denote this disk with smooth boundary by $\Delta$.  So,
$\partial\Delta\subset\tau$ and the directions of the arcs in
$\partial\Delta$ either all point inwards or all point outwards, as
$\tau$ is transversely oriented according to the branch direction.
Since we have assumed that $D'$ contains no subdisk that is a
generalized sink disk, the direction of $\partial\Delta$ must point
out of $\Delta$, and hence $\Delta\subset int(D')$.  Therefore, by our
construction of $G$, $\Delta$ must be parallel to a disk component of
$\partial_hN(B)$ (see Definition~\ref{D:parallel} for the definition
of parallel).  After an isotopy in a small neighborhood of $\Delta$,
we can assume that $\Delta$ is a disk component of $\partial_hN(B)$.
Since $\partial_hN(B)$ is incompressible, $\Delta$ must be a
horizontal boundary component of a $D^2\times I$ region of
$M-int(N(B))$.  Let $K=D^2\times I$ ($I=[-1,1]$) be the component of
$M-int(N(B))$ such that $\Delta=D^2\times\{-1\}\subset\partial K$.  We
denote $D^2\times\{1\}\subset\partial K$ by $\Delta'$.  Then, we can
isotope $D'$ across $K$ in a small neighborhood of $K$.  In other
words, $(D'-\Delta)\cup A\cup\Delta'$, where $A=\partial D^2\times
I\subset\partial K$ is the vertical boundary of $K$, is an embedded
disk in $N(B)$ that is isotopic (in $M$) to $D'$.  Then, by a small
perturbation near $A$, we can isotope the disk $(D'-\Delta)\cup
A\cup\Delta'$ to be transverse to the $I$--fibers of $N(B)$.  We denote
the disk after this perturbation by $D''$.  Clearly $D''$ is isotopic
to $D'$.  Moreover, we can assume that $\Delta'\subset D''$ and $D'$
coincides with $D''$ outside a small neighborhood of $\Delta$.  The
picture of $D'\cup D''$ is like a disk with an ``air bubble" inside
which corresponds to the $D^2\times I$ region $K$.

\begin{figure}[ht!]
\centerline{\small
\SetLabels 
\E(0.1*0.5){$G$}\\
\E(0.95*0.5){$\tau$}\\
\endSetLabels 
\AffixLabels{{\includegraphics[width=4in]{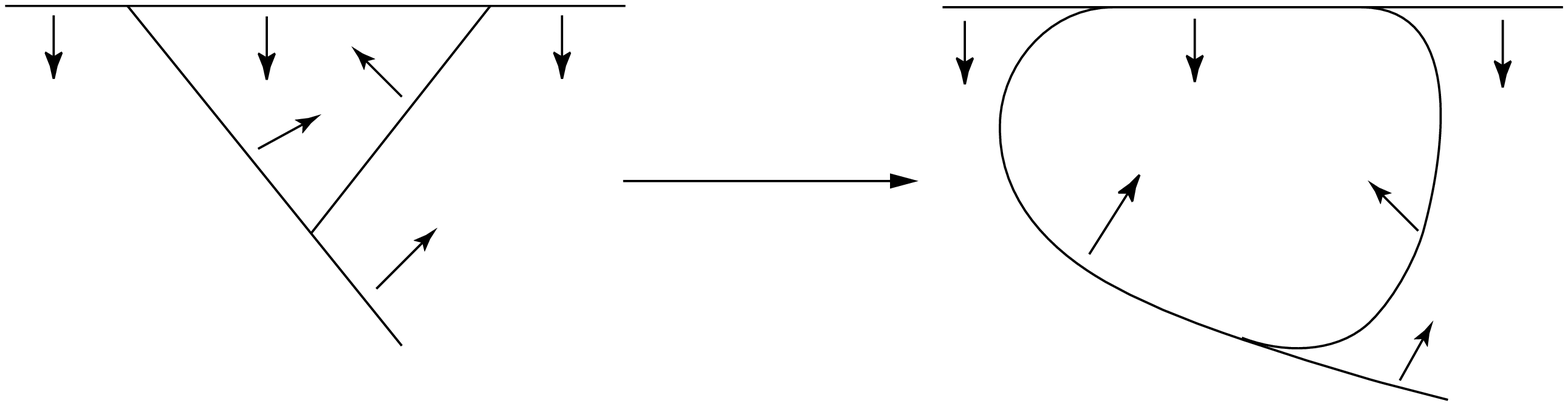}}}}
\vspace{2mm}
\caption{}\label{traintrack}
\end{figure}

Let $G'$ be the union of $x\in int(\Delta)$ such that
$J_x\subset\pi^{-1}(L)$ and $\partial_vN(B)\cap J_x$ contains a
component on the negative side of $x$.  Then, $G'\cup\partial\Delta$
is also a trivalent graph, and each edge of $G'$ has a direction
induced from the branch direction of the corresponding arc in $L$.
Note that $G'\cup\partial\Delta=\Delta\cap\pi^{-1}(L)$, since
$\Delta\subset\partial_hN(B)$.  As before, we can deform
$G'\cup\partial\Delta$ into a transversely oriented train track
$\tau'$.  By a standard index argument, $\Delta-\tau'$ must have a
smooth disk component, ie, there is a smooth circle in $\tau'$ which
is the boundary of the closure of a component of $\Delta-\tau'$.  We
denote this disk with smooth boundary by $\delta$.  Since $\tau'$ is
transversely oriented and $\partial\delta$ is a smooth circle in
$\tau'$, the directions of the arcs in $\partial\delta$ either all
point into $\delta$ or all point out of $\delta$.  If the direction of
$\partial\delta$ points into $\delta$, since
$G'\cup\partial\Delta=\Delta\cap\pi^{-1}(L)$, $\pi(\delta)$ is a sink
disk, which gives a contradiction.  Thus, $\pi(\delta)$ must be a
branch of $B$ with branch direction of each boundary arc pointing
outwards.  Moreover, since $\pi|_{D'}$ is injective, $\pi(\delta)$ is
a removable disk.

Next, we show that $\pi(D'')$ does not contain $\pi(int(\delta))$, and
hence $D''$ is carried by the branched surface $B-int(\delta)$, as
$\delta$ is removable.  We first show that there is no $I$--fiber of
$N(B)$ that intersects both $\Delta'$ and $int(\delta)$ (note that
$\Delta'\subset D''$ and $\delta\subset\Delta\subset D'$).  Otherwise,
since $\delta\cap\pi^{-1}(L)=\partial\delta$, $\delta$ is parallel to
a subdisk of $\Delta'$.  As $\partial\delta\subset\partial\Delta\cup
G'$, we have two cases: one case is $\partial\delta=\partial\Delta$
and the other case is $\partial\delta\cap G'\cap
int(\Delta)\ne\emptyset$.  If $\partial\delta=\partial\Delta$, we have
$\delta=\Delta$ and $G'=\emptyset$.  Then, since $\delta$ is parallel
to a subdisk of $\Delta'$ and $G'=\emptyset$, $\Delta=\delta$ is
parallel to a subdisk in the interior of $\Delta'$.  Since $\Delta$
and $\Delta'$ are the two components of the horizontal boundary of the
$D^2\times I$ region $K$, $K$ forms a standard Reeb component, which
gives a contradiction.  In fact, in this case, $B-\pi(int(\delta))$ is
a branched surface with one horizontal boundary component a torus that
bounds a solid torus in $M$.  Thus, $\partial\delta\cap G'\cap
int(\Delta)\ne\emptyset$. As $\delta$ is parallel to a subdisk of
$\Delta'$ and $\delta\subset\Delta$, there is an $I$--fiber $J$ of
$N(B)$ that intersects both $\Delta'$ and $\partial\delta\cap
int(\Delta)$.  Note that since $\Delta$ is a disk component of
$\partial_hN(B)$ and $J\cap int(\Delta)\ne\emptyset$, one endpoint of
$J$ must lie in $int(\Delta)$.  As $\partial\delta\subset\pi^{-1}(L)$,
$J\cap\partial_vN(B)\ne\emptyset$.  Moreover, since $\Delta'$ is also
a disk component of $\partial_hN(B)$ and since $\delta$ is parallel to
a subdisk of $\Delta'$, $J\cap\partial\Delta'\ne\emptyset$.  Then, as
$\Delta$ and $\Delta'$ are the two horizontal boundary components of
the $D^2\times I$ region $K$, $J\cap\partial\Delta'\ne\emptyset$
implies $J\cap\partial\Delta\ne\emptyset$.  Thus, $J\cap
int(\Delta)\ne\emptyset$ and $J\cap\partial\Delta\ne\emptyset$, and
hence $J\cap\Delta$ contains at least two points, which contradicts
our assumption that $\pi|_{D'}$ is injective.  Therefore, there is no
$I$--fiber of $N(B)$ that intersects both $\Delta'$ and $int(\delta)$.

Since $\pi|_{D'}$ is injective and since there is no $I$--fiber of
$N(B)$ that intersects both $\Delta'$ and $int(\delta)$, by our
construction of $D''$, there is no $I$--fiber of $N(B)$ that intersects
both $D''$ and $int(\delta)$.  Since $\pi(\delta)$ is a removable
disk, $D''$ is a generalized sink disk in $N(B^-)$, where $B^-$ is the
branched surface $B-\pi(int(\delta))$.  Note that $\pi|_{D''}$ is not
necessarily injective.

Now, $D''$ is a generalized sink disk in a fibered neighborhood of the
branched surface $B^-=B-\pi(int(\delta))$.  We can then apply the same
argument above to $B^-=B-\pi(int(\delta))$, replacing $B$ and $D$ by
$B^-$ and $D''$ respectively.  As in the argument above, the existence
of a generalized sink disk always yields a removable disk (such as the
$\delta$ above).  However, if we keep eliminating these removable
disks, we eventually get an efficient laminar branched surface that
still has a generalized sink disk.  This gives a contradiction.
\end{proof}

\begin{remark}
It is easy to see from the proof of Lemma~\ref{L:gsd} that if a
branched surface $B$ contains a trivial bubble but has no sink disk,
then $B$ must contain a removable disk.
\end{remark}

An easy corollary (Corollary~\ref{C:contact}) of Lemma~\ref{L:gsd} is
that there is no disk of contact in a laminar branched surface.
Figure~\ref{F13} (a) is the simplest example of a disk of contact.  By
definition (see condition 4 in Proposition~\ref{P11}), a disk of
contact is an embedded disk $D\subset N(B)$ such that $D$ is
transverse to the $I$--fibers of $N(B)$ and $\partial
D\subset\partial_v N(B)$.  If $\pi^{-1}(L)\cap int(D)\ne\emptyset$,
$\pi(D)$ is not even a branch of $B$.  For example, we can add some
branches to Figure~\ref{F13} (a) in a complicated way, but it can
still be a disk of contact by definition.  In general, it is not
obvious that the condition of no sink disks implies that there is no
disk of contact, although Figure~\ref{F13} (a) is an example of sink
disk.

\begin{corollary}\label{C:contact}
A laminar branched surface does not contain any disk of contact.
\end{corollary}
\begin{proof}
By the definition of disk of contact in Proposition~\ref{P11}, a disk
of contact is a generalized sink disk and Corollary~\ref{C:contact}
follows from Lemma~\ref{L:gsd}.
\end{proof}

The next two lemmas will be used in section~\ref{S:construction}, and
the proofs are essentially the same as the proof of Lemma~\ref{L:gsd}.

\begin{lemma}\label{L:dl}
Let $B$ be a laminar branched surface.  Then, $N(B)$ contains no disk
$D$ with the following properties:
\begin{enumerate}
\item $D$ is an embedded disk in $N(B)$ that is transverse to the
$I$--fibers of $N(B)$;
\item $\pi(\partial D)$ is a nontrivial simple closed curve in $B-L$.
\end{enumerate}
\end{lemma}
\begin{proof}
We first show that if $N(B)$ contains such a disk $D$, then $D$ has a
subdisk $E$ such that $\pi(\partial E)=\pi(\partial D)$ and
$\pi(\partial E)\cap\pi(int(E))=\emptyset$.  Since $\pi(\partial D)$
is a nontrivial simple closed curve in $B-L$,
$D\cap\pi^{-1}(\pi(\partial D))$ is a union of simple closed curves in
$D$.  Let $E\subset D$ be a disk bounded by an innermost (among curves
in $D\cap\pi^{-1}(\pi(\partial D))$) simple closed curve.  Then, since
$\partial E$ is innermost, $\pi(\partial E)\cap\pi(int(E))=\emptyset$.
Therefore, we may assume that our disk $D$ has an additional property
that $\pi(\partial D)\cap\pi(int(D))=\emptyset$.

If $\pi|_{D}$ is not injective, since $\pi(\partial
D)\cap\pi(int(D))=\emptyset$, similar to the proof of
Lemma~\ref{L:gsd}, there must be a subdisk in $int(D)$ that is a
generalized sink disk, which contradicts Lemma~\ref{L:gsd}.  More
precisely, let $n$ be the maximal number of intersection points of $D$
with any $I$--fibers of $N(B)$, and $X_n$ be the union of $I$--fibers of
$N(B)$ whose intersection with $D$ consists of $n$ points.  Since
$\pi|_{D}$ is not injective, $n>1$.  Then, since $\pi(\partial
D)\cap\pi(int(D))=\emptyset$ and $\pi|_{\partial D}$ is injective,
$X_n\cap D$ is a collection of compact subsurfaces of $int(D)$.
Moreover, since $n$ is maximal, the boundary of $X_n\cap D$ lies in
$\pi^{-1}(L)$ with direction (induced from the branch direction of
$L$) pointing into $X_n\cap D$.  Thus, the outer boundary of a
component of $X_n\cap D$ bounds a generalized sink disk in $int(D)$,
which contradicts Lemma~\ref{L:gsd}.  Therefore, $\pi|_{D}$ must be
injective.

Since $\pi|_{D}$ is injective, as in the proof of Lemma~\ref{L:gsd},
we can find a removable disk $\delta$ in $int(D)$.  Moreover, we can
find another disk $D'$, which we get by isotoping $D$ across a
$D^2\times I$ region, such that $\partial D=\partial D'$ and
$\pi(D')\cap\pi(int(\delta))=\emptyset$.

Therefore, $D'$ satisfies the two hypotheses (for $D$) in the lemma,
and $D'$ is carried by the branched surface $B-int(\delta)$.  Then, we
can apply the same argument to the branched surface $B-int(\delta)$,
replacing $B$ and $D$ by $B-int(\delta)$ and $D'$ respectively.
Similar to the proof of Lemma~\ref{L:gsd}, we get a contradiction once
the branched surface becomes efficient.
\end{proof}

\begin{lemma}\label{L:ba}
Let $B$ be a laminar branched surface.  Then, $N(B)$ contains no disk
$D$ with the following properties:
\begin{enumerate}
\item $D$ is an embedded disk in $N(B)$ that is transverse to the
$I$--fibers of $N(B)$;
\item $\pi(\partial D)$ is a simple closed curve in $B$ that is
transverse to $L$ ($L\cap\pi(\partial D)\ne\emptyset$) and does not
contain any double point of $L$;
\item the points in $L\cap\pi(\partial D)$ have coherent branch
directions along $\pi(\partial D)$ (clockwise or counterclockwise),
where we consider the branch direction of each point in
$L\cap\pi(\partial D)$ to be along $\pi(\partial D)$, ie, a small
neighborhood of $\pi(\partial D)$ is either a branched annulus or a
branched M\"{o}bius band with coherent branch direction as shown in
Figure~\ref{F34}.
\end{enumerate}
\end{lemma}
\begin{proof}
The proof is very similar to the proof of Lemma~\ref{L:dl}.  We first
show that $D\cap\pi^{-1}(\pi(\partial D))$ is a union of simple closed
curves in $D$.  Since $\pi^{-1}(\pi(\partial D))$ is a compact set,
$D\cap\pi^{-1}(\pi(\partial D))$ is a union of circles or compact arcs
in $D$.  If $D\cap\pi^{-1}(\pi(\partial D))$ has a component that is a
compact arc, which we denote by $\alpha$, then by our hypothesis that
$\pi|_{\partial D}$ is injective, $\partial\alpha$ must lie in
$\pi^{-1}(L)$ with direction (consistent with the branch direction)
pointing into $\alpha$.  However, this is impossible because the
points in $L\cap\pi(\partial D)$ have coherent branch directions along
$\pi(\partial D)$, in other words, there is no subarc of $\pi(\partial
D)$ with endpoints in $L$ and branch directions of both endpoints
pointing into this arc.  Thus, $D\cap\pi^{-1}(\pi(\partial D))$ is a
union of simple closed curves in $D$.  Let $E\subset D$ be a disk
bounded by an innermost (among curves in $D\cap\pi^{-1}(\pi(\partial
D))$) simple closed curve.  Since $\partial E$ is innermost,
$\pi(\partial E)\cap\pi(int(E))=\emptyset$.  Therefore, similar to the
proof of Lemma~\ref{L:dl}, we can assume that our disk $D$ has an
additional property that is $\pi(\partial
D)\cap\pi(int(D))=\emptyset$.

Thus, as in the proof of Lemma~\ref{L:dl}, $\pi|_{D}$ must be
injective.  Then, as in the proof of Lemma~\ref{L:gsd}, we can
construct a train track $\tau\subset int(D)$ as follows.  We first fix
a normal direction for $D$.  For every point $x\in int(D)$, let $J_x$
be the $I$--fiber of $N(B)$ that contains $x$.  Then, $J_x-x$ has two
components.  According to the fixed normal direction of $D$, we say
that the points in one component of $J_x-x$ are on the positive side
of $x$, and points in the other component of $J_x-x$ are on the
negative side of $x$.  Let $G$ be the union of $x\in int(D)$ such that
$J_x\subset\pi^{-1}(L)$ and $\partial_vN(B)\cap J_x$ contains a
component on the positive side of $x$.  Then, $G$ is a trivalent graph
and each edge has a direction induced from the branch direction.  As
shown in Figure \ref{traintrack}, this trivalent graph $G$ can be
deformed into a transversely oriented train track $\tau$ according to
the directions of the edges in $G$.  By fixing an appropriate normal
direction for $D$, we can assume that $\tau\ne\emptyset$.

Since the branch directions of points in $L\cap\pi(\partial D)$ are
coherent along $\pi(\partial D)$ and since $\tau$ is transversely
oriented according to the branch direction, there is no arc carried by
$\tau$ with both endpoints in $\partial D$.  Then, similar to the
argument in the Poincar\'{e}-Bendixson theorem \cite{N}, $\tau$ must
carry a circle that bounds a disk in $int(D)$.  Hence, there must be a
smooth disk whose boundary is a smooth circle in $\tau$ and whose
interior is a component of $int(D)-\tau$.  As in the proof of
Lemma~\ref{L:gsd}, we can find a removable disk in $int(D)$, and we
can get another disk $D'$ (with $\partial D=\partial D'$) by isotoping
$D$ across a $D^2\times I$ region $K$ of $M-int(N(B))$.  Moreover,
$\pi(D')$ does not pass through the removable disk $\delta$.

Then, we can apply the argument above again to the branched surface
$B-int(\delta)$, replacing $B$ and $D$ by $B-int(\delta)$ and $D'$
respectively.  As in the proof of Lemma~\ref{L:gsd}, we get a
contradiction once the branched surface becomes efficient.
\end{proof}

\section{Extending laminations}\label{S:extend}

In this section, we show that, in most cases, we can extend a
lamination from the vertical boundary of an $I$--bundle over a surface
to its interior.  The results in this section appear in \cite{G4}
implicitly, and most of the proof we give here is in fact a
modification of the arguments in \cite{G4}.

Let $B$ be a branched surface carrying a lamination $\lambda$.
Suppose $\partial B$ is a union of circles.  By `blowing air' into
leaves, ie, replacing leaves by $I$--bundles over these leaves and
deleting the interior of these $I$--bundles, we can assume that
$\lambda$ is nowhere dense in $N(B)$.  For simplicity, we will assume
the intersection of $\lambda$ with every interval fiber is a Cantor
set.

Let $I=[-1, 1]$ and $Homeo^+(I)$ be the group of self-homeomorphism of
$I$ fixing endpoints.  The next lemma is well-known, and the proof is
easy (see also \cite{EHN}).

\begin{lemma}\label{L21}
Any map $f\in Homeo^+(I)$ is a commutator, ie, there are $g, h\in
Homeo^+(I)$ such that $f=g\circ h\circ g^{-1}\circ h^{-1}$.
\end{lemma}
\begin{proof}
As the fixed points of $f$ is a closed set in $I$ and the complement
of a closed set is a union of intervals, it suffices to prove
Lemma~\ref{L21} for maps without fixed points in the interior of $I$.
Hence, we may assume that $f(z)>z$ for any $z\in (-1,1)$ (the case
$f(z)<z$ is similar).  It suffices to show that $f$ is conjugate to
any map $p\in Homeo^+(I)$ with the property that $p(z)>z$ for any
$z\in (-1,1)$.

Let $x$ be an arbitrary point in the interior of $I$.  As $f(x)>x$ and
$p(x)>x$, the intervals $[f^n(x), f^{n+1}(x)]$ ($n\in\mathbb{Z}$,
$f^0(x)=x$ and $f^1=f$) partition the interval $I$, and the intervals
$[p^n(x), p^{n+1}(x)]$ ($n\in\mathbb{Z}$) also partition the interval
$I$.  Let $q_0\co  [x, f(x)]\to [x, p(x)]$ be any homeomorphism fixing
endpoints.  We define $q_n=p^n\circ q_0\circ f^{-n}\co  [f^n(x),
f^{n+1}(x)]\to [p^n(x), p^{n+1}(x)]$.  These maps $q_n$ fit together
to give a homeomorphism $q\co [-1, 1]\to [-1,1]$, and it follows from the
definition of $q_n$ that $f=q^{-1}\circ p\circ q$, ie, $f$ and $p$
are conjugate.
\end{proof}

The following lemma is an application of Lemma~\ref{L21}.

\begin{lemma}\label{L22}
Let $c$ be a circular component of $\partial B$.  If $B$ fully carries
a lamination, then the new branched surface constructed by gluing $B$
and a once-punctured orientable surface with positive genus along $c$
also carries a lamination.
\end{lemma}

\begin{proof}
Let $A=\pi^{-1}(c)$, where $\pi \co N(B)\to B$ is the collapsing map.
Then $\lambda |_A$ is a one-dimensional lamination in the annulus
$int(A)=S^1\times int(I)$, and $A-\lambda |_A$ is a union of
$I$--bundles.  Each $I$--bundle is homeomorphic to either
$\mathbb{R}\times I$ or $S^1\times I$.  We trivially extend $\lambda
|_A$ to a (one-dimensional) foliation of $int(A)$ by associating each
$I$--bundle with its canonical product foliation.  Assume that the
foliation of $int(A)$ constructed above is the suspension of a
homeomorphism $f\co  int(I)\to int(I)$

Let $S$ be the once-punctured surface that we glue to $B$.  We
consider the $I$--bundle $S\times I$ and $A=\partial S\times I=c\times
I$.  By Lemma 2.1, there exist $a_1, b_1,\dots, a_g, b_g$ such that
$f=[a_1,b_1]\circ\dots\circ [a_g,b_g]$, where $a_i, b_i$ are
homeomorphisms of $int(I)$ and $g$ is the genus of $S$.  By attaching
thick bands foliated by the suspensions of $a_i$'s and $b_i$'s to a
$disk\times I$ with the trivial product foliation, we can build
$S\times I$.  The foliations of the thick band and the disk can be
glued together according to the identity map of $I$.  This gives us a
foliation of $S\times I$ whose boundary on $\partial S\times I$ is the
suspension of $f$.  In other words, we can extend the foliation of $A$
to a foliation of $S\times I$.

Then by `blowing air' into leaves, we can change the foliation of
$S\times I$ to a nowhere dense lamination $\nu$ such that $\lambda
|_A$ is a sub-lamination of $\nu |_A$.  Indeed, by our construction of
the foliation of $A$, $\nu |_A$ is just $\lambda |_A$ plus some
parallel nearby leaves.  Now we change the lamination $\lambda$ in
$N(B)$ by adding some parallel leaves so that the new lamination
restricted to $A$ is the same as $\nu |_A$.  Gluing up the two
laminations, we get a lamination fully carried by the new branched
surface.
\end{proof}

\begin{remark}
The operations we used on laminations and foliations in the proof
above are standard, see operations 2.1.1, 2.1.2, 2.1.3 in \cite{G4}.
\end{remark}

\begin{corollary}\label{C23}
Let $c_1, c_2,\dots, c_n$ be n circular components of $\partial B$.
If $B$ fully carries a lamination, then the new branched surface
constructed by gluing a non-planar orientable surface with n boundary
components along $c_i$'s fully carries a lamination.
\end{corollary}

\begin{proof}
We first glue a planar surface with $n+1$ boundary components to $B$.
By adding thickened bands between $c_1, c_2,\dots, c_n$, we can
trivially extend the lamination through the planar surface.  Then we
can glue a once-punctured surface to the $(n+1)$th boundary component
of the planar surface and the result follows from Lemma~\ref{L22}.
\end{proof}

The next Lemma is a modification of operation 2.4.4 in \cite{G4}.

\begin{lemma}\label{L24}
Let $c_1$ and $c_2$ be two circular components of $\partial B$.  If
$B$ fully carries a lamination without disk leaves, then the new
branched surface constructed by gluing an annulus between $c_1$ and
$c_2$ carries a lamination.
\end{lemma}

\begin{proof}
Let the vertical boundary components of $N(B)$ along $c_1$ and $c_2$
be annuli $A_1$ and $A_2$, $A_i=c_i\times [-1,1]$.  What we want to do
is to add some leaves to $\lambda$ so that the restriction of the new
lamination to $A_1$ and $A_2$ are the same, hence we can glue them
together.

First, we replace every boundary leaf of $\lambda$ by an embedded
$I$--bundle over this leaf, then we delete the interior of the
$I$--bundle.  We still call this lamination $\lambda$.  After this
operation, $\lambda |_{A_i}$ has two pairs of isolated circles near
$\partial A_i$.

Then we isotope $\lambda$ such that $\lambda |_{c_1\times [-1,0]}$ is
an isolated circle, say $e_1$, $e_1=c_1\times \{-1\}$, and $\lambda
|_{c_2\times [0,1]}$ is also an isolated circle, say $e_2$,
$e_2=c_2\times\{1\}$.  Let $L_1$ and $L_2$ be the leaves in $\lambda$
corresponding to $e_1$ and $e_2$ respectively.  Clearly $L_1$ and
$L_2$ are orientable surfaces.  Then we add two leaves $L_1', L_2'$ to
$\lambda$ which are parallel and close to $L_1, L_2$ respectively such
that $L_i\cup L_i'$ bounds a product region in $N(B)$.  This is
actually the same operation as replacing $L_i$ by an $I$--bundle and
deleting the interior of the $I$--bundle.  Let $L_i'\cap A_i=e_i'$,
$i=1,2$.  Let the annulus in $A_i$ bounded by $e_1'\cup c_1\times
\{1\}$ be $J_1$, the one bounded by $c_2\times\{-1\}\cup e_2'$ be
$J_2$, the one bounded by $e_1'\cup e_1$ be $K_1$, and the one bounded
by $e_2'\cup e_2$ be $K_2$.

Before we proceed, we point out a fact that is the following.  Let
$F\times I$ be a product region over a surface $F$ ($F$ could be
non-compact).  Suppose $F$ is not a disk and $C$ is a boundary
component of $F$.  Then any foliation on $C\times I$ can be extended
to the whole of $F\times I$.  The proof is easy.  If $F$ has another
boundary component or an end, the construction is trivial, and
otherwise, it follows from Lemma~\ref{L21}.

\medskip
\noindent
\textbf{Case 1}\qua One of $L_1$ and $L_2$ (say $L_1$) is not a compact
planar surface with boundary on $A_1\cup A_2$.

\medskip
\noindent
\textbf{Case 1a}\qua  $L_2$ is not a compact planar surface with boundary
on $A_1\cup A_2$ either.

We foliate $J_1$ and $J_2$ as before, ie, foliate all annular
components of $A_i-\lambda |_{A_i}$ by circles and other components by
adding spirals coherent to $\lambda |_{A_i}$.  Then we foliate $K_1$
with the same foliation of $J_2$ and $K_2$ with the same foliation of
$J_1$.  Now the foliation on $A_1$ and $A_2$ are the same.  By our
assumption on $L_1$ and $L_2$, we can extend the foliation of $K_i$ to
the product region bounded by $L_i\cup L_i'$.  Then, as before, by
`blowing air' into leaves we can change the foliation on $A_i$ to be a
nowhere dense lamination that contains $\lambda |_{A_i}$ as a
sub-lamination.  By our construction of foliation on $A_i$, the
complement of $\lambda |_{A_i}$ is a product lamination.  After
possibly replacing every leaf by a product lamination of $leaf\times\{
a\ cantor\ set\}$, we can extend the lamination on $A_i$ to $N(B)$.
Now the new lamination in $N(B)$ when restricted to $A_1$ and $A_2$,
gives the same lamination.

\medskip
\noindent
\textbf{Case 1b}\qua $L_2$ is a compact planar surface with boundary on
$A_1\cup A_2$, but $L_2\cap J_1=\emptyset$.

This case is very similar to Case 1a.  We first foliate $J_1$ in the
same way as before, then give $K_2$ the same foliation as that of
$J_1$.  Since $L_2$ is not a disk, we can extend the foliation on
$K_2$ to the product region bounded by $L_2\cup L_2'$.  Now we might
have changed the lamination on $J_2$.  We can extend the (new)
lamination on $J_2$ to a foliation as before and give $K_1$ the same
foliation as $J_2$, and the rest is as in Case 1a.

Before we proceed, we quote the Lemma 2.1 of \cite{G4}.

\begin{lemma}\label{L25}
Let $f$, $h$, $\sigma$, $\tau$ be either homeomorphisms of $I$ fixing
endpoints or maps of the empty set.

{\rm i)}\qua There exists a homeomorphism $g$ conjugate to the concatenation of
$f$, $g$, $h$.

{\rm ii)}\qua There exists homeomorphisms $g$, $\mu$ of $I$ such that $\mu$ is
conjugate to the concatenation of $f$, $g^{-1}$and $h$, and $g$ is
conjugate to the concatenation of $\sigma$, $\mu^{-1}$ and $\tau$.
\end{lemma}

\begin{remark}
Let $A$ be an annulus and $\mathcal{F}_1, \mathcal{F}_2$ be two
foliations on $\partial A\times I$.  Suppose $\mathcal{F}_i$ is a
suspension of a homeomorphism of $I$ fixing endpoints, say $f_i$,
$i=1,2$.  Then we can extend $\mathcal{F}_1$ and $\mathcal{F}_2$ to a
foliation of $A\times I$ if and only if $f_1$ is conjugate to $f_2$.
\end{remark}

\noindent
\textbf{Case 1c}\qua $L_2$ is a compact planar surface and has some
boundary component $E$ in $J_1$.

Since $L_2'\cup L_2$ bounds a product region, $L_2'$ has a boundary
component $E'$ and $E'\cup E$ bounds an annulus $J'$ in $J_1$.  We
first extend the lamination on $J_1-J'$ to a foliation as before, and
assume that this foliation is a suspension of maps $f$ and $h$ (since
$J_1-J'$ consists of two annuli).  Then we construct the same
foliation, which is the suspension of $g$, on $J'$ and $K_2$, where
$g$ is as in Lemma~\ref{L25}~ (i).  By Lemma~\ref{L25}~ (i), $K_2$ and
$J_1$ have the same foliation.  So we can extend it to a foliation in
the product region bounded by $L_2\cup L_2'$, and the rest is the same
as Case 1b.

\medskip
\noindent
\textbf{Case 2a}\qua  Both $L_1$ and $L_2$ are planar and some non-$e_1$
component of $\partial L_1$ is disjoint from $J_1$.

We first foliate $J_1$ as before, then give $K_2$ the same foliation
as that of $J_1$ and extend it to a foliation of the product region
bounded by $L_2\cup L_2'$.  Applying Lemma~\ref{L25}~(i), if
necessary, we can construct the same foliation on $K_1$ and $J_2$ such
that it can be extended to a foliation in the product region bounded
by $L_1\cup L_1'$ (using our assumption of $\partial L_1$).

\medskip
\noindent
\textbf{Case 2b}\qua Both $L_1$ and $L_2$ are compact planar surfaces and
all the non-$e_i$ components of $\partial L_i$ are in $J_i$, $i=1,2$.

Let $d_i$ be another boundary component of $L_i$ and $d_i'$ be the
corresponding boundary component of $L_i'$.  Then $d_i\cup d_i'$
bounds a annulus $J_i'$ in $J_i$, $i=1,2$.  We extend the lamination
on $J_i-J_i'$ as before, and foliate $K_1$ by the suspension of a map
$g$, $J_1'$ by the suspension of map $g^{-1}$, $K_2$ by the suspension
of a map $\mu$, and $J_2'$ by the suspension of map $\mu^{-1}$.  By
our assumption on $L_1$ and $L_2$, we can extend the foliation to the
product region bounded by $L_i\cup L_i'$, $i=1,2$.  Using
Lemma~\ref{L25}~(ii), we can find maps $g$ and $\mu$ such that the
foliation on $J_i$ is the same as the foliation on $K_j$, $i\ne j$,
and the rest is the same as before.
\end{proof}

\begin{lemma}\label{L26}
Let $c_1, c_2,\dots, c_n$ be n circular components of $\partial B$.
If $B$ fully carries a lamination without disk leaves, then the new
branched surface constructed by gluing a non-disk surface with n
boundary components along $c_i$'s fully carries a lamination.
\end{lemma}
\begin{proof}
Let $S$ be the surface that we glue to $B$.  The case that $S$ is
orientable follows from the lemmas above.  As in the previous
arguments, we only need to consider case that $S$ is a M\"{o}bius
band.

Let $c_1=\partial S$ and $A=c_1\times [-1, 1]$.  As in the proof of
Lemma~\ref{L24}, by replacing a boundary leaf by an $I$--bundle over
this leaf and deleting the interior of the $I$--bundle, we can assume
that $c_1\times\{-1, 0, 1\}\subset\lambda |_A$ are isolated circles.
Since the ambient manifold $M$ is assumed to be orientable, we can
glue a twisted $I$--bundle over a M\"{o}bius band, which we denote by
$U$, to $N(B)$ along $A=c_1\times [-1,1]$.  Then, $c_1\times\{ 0\}$
bounds a M\"{o}bius band $u$ in $U$.  Topologically, $U-u$ is a
fibered neighborhood of an annulus with each fiber a half open and
half closed interval, ie, $U-u=annulus\times [a,b)$.  The vertical
boundary of $U-u$ is the union of $c_1\times [-1,0)$ and $c_1\times
(0,1]$.  By Lemma~\ref{L24}, we can extend the lamination through
$U-u$ and the Lemma holds.
\end{proof}

\section{Constructing laminations carried by branched surfaces}\label{S:construction}

\begin{figure}[ht!]
\centerline{\small
\SetLabels 
\E(0.36*0.1){$L'$}\\
\E(0.6*0.15){$P(L')$}\\
\E(0.15*-0.1){(a)}\\
\E(.7*-0.1){(b)}\\
\endSetLabels 
\AffixLabels{{\includegraphics[width=4in]{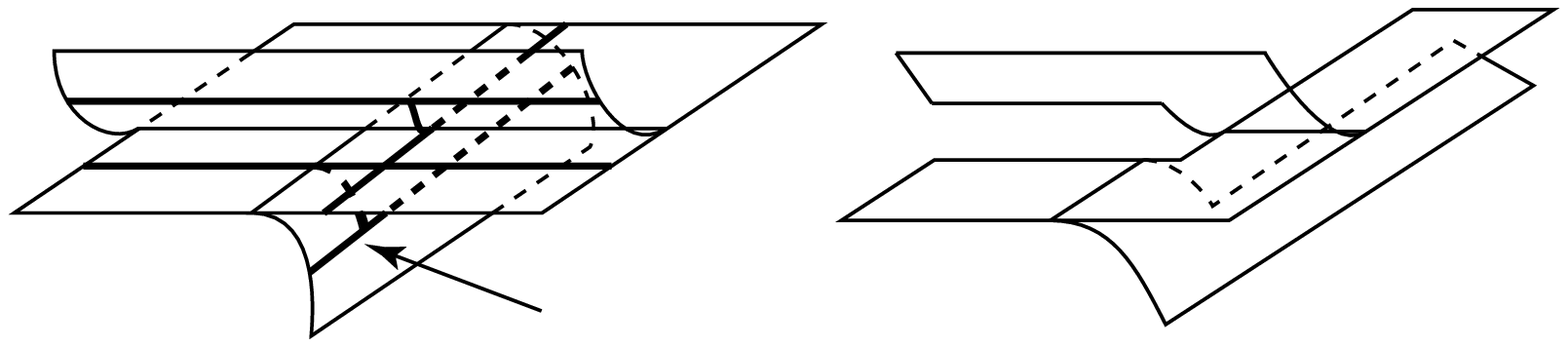}}}}
\vspace{2mm}
\caption{}\label{F31}
\end{figure}

Suppose $B$ is a laminar branched surface.  Let $L'$ be a graph in
$B$, whose local picture is as shown in Figure~\ref{F31} (a).  We can
also describe $L'$ as follows.  Let $l_1, l_2,\dots,l_s$ be the
boundary curves of the surface $\partial_hN(B)$.  For each $l_i$, we
take a simple closed curve $l_i'$ in the interior of $\partial_hN(B)$
that is close and parallel to $l_i$.  Let $DL=\cup_{i=1}^s\pi (l_i)$.
Near every double point of $L$, the intersection of $DL$ with $L$
consists of two points.  Then, we add some short arcs connecting these
intersection points to $DL$, as shown in Figure~\ref{F31} (a), and the
union of $DL$ and these short arcs is $L'$.

Let $K_{L'}$ be a closed small regular neighborhood of $L'$ in $M$.
Let $P(L')=B\cap K_{L'}$, whose local picture is as shown in
Figure~\ref{F31} (b).  We call $B\cap\partial K_{L'}$ the
\emph{boundary} of the branched surface $P(L')$, and denote $B\cap
int(K_{L'})$ by $int(P(L'))$.  The branch locus of the branched
surface $P(L')$ is a union of simple arcs, as shown in
Figure~\ref{F31} (b).

\begin{figure}[ht!]
\centerline{\small
\SetLabels 
\E(0.23*.52){$D$}\\
\E(0.75*-0.1){$D^B$}\\
\E(0.2*-0.3){(a)}\\
\E(.8*-0.3){(b)}\\
\endSetLabels 
\AffixLabels{{\includegraphics[width=4in]{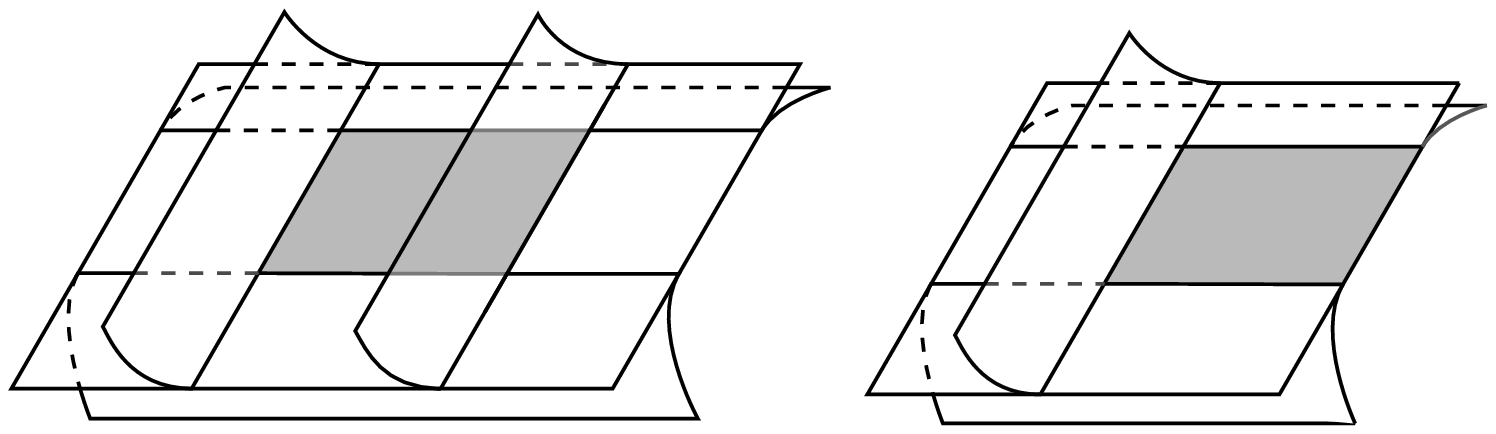}}}}
\vspace{6mm}
\caption{}\label{F32}
\end{figure}

There is a one-one correspondence between the components of $B-L$ and
the components of $B-P(L')$.  For each branch $D$ of $B$, we denote
the corresponding component of $B-int(P(L'))$ by $D^B$.  For example,
the shaded branch in Figure~\ref{F32} (a) corresponds to
Figure~\ref{F32} (b), which is a component of $B-int(P(L'))$.  The
relation between $D$ and $D^B$ can also be described as follows.  We
consider $N(B)$ as a fibered regular neighborhood of $B$, and $B$ lies
in the interior of $N(B)$ such that every $I$--fiber of $N(B)$ is
transverse to $B$.  To simplify notation, we do not distinguish the
$B$ in the interior of $N(B)$ and the $B$ as the image of the map
$\pi\co N(B)\to B$, which collapses every $I$--fiber of $N(B)$ to a point.
There is a natural one-one correspondence between components of $B-L$
and components of $N(B)-\pi^{-1}(L)$.  For any branch $D$ of $B$,
$int(D)$ is a component of $B-L$, and the corresponding component of
$N(B)-\pi^{-1}(L)$ is an $I$--bundle over $int(D)$, whose intersection
with the $B$ lying in $int(N(B))$ is the same as the component of
$B-P(L')$ that corresponds to $int(D)$.

For any branch $D$ of $B$, we denote by $N_B(D)$ the closure in the
path metric of the component of $N(B)-\pi^{-1}(L)$ that corresponds to
$D$.  Thus, $N_B(D)$ is an $I$--bundle over $D$ with a bundle structure
induced from $N(B)$.  The vertical boundary of $N_B(D)$ is bundle
isomorphic to $\partial D\times I$, and we can identify
$N_B(D)-\partial D\times I$ with the component of $N(B)-\pi^{-1}(L)$
that corresponds to $D$.  By our argument above, $B\cap
(N_B(D)-\partial D\times I)$ is the same as the component of $B-P(L')$
corresponding to $D$.  Thus, we can assume that $D^B$, the
corresponding component of $B-int(P(L'))$ lies in $N_B(D)$ with
$\partial D^B\subset\partial D\times I$.

We can reconstruct $N(B)$ by gluing all the $N_B(D)$'s (for all the
branches of $B$) together along their vertical boundaries, and
simultaneously, those $D^B$'s (lying in $N_B(D)$'s) are glued together
to form $B$.  Moreover, the gluing (for $D^B$'s) above is essentially
the same as gluing the $D^B$'s and $P(L')$ together to form $B$.

\begin{figure}[ht!]
\centerline{\small
\SetLabels 
\E(0.5*0.73){$\{p\}\times I$}\\
\E(0.5*0.16){$E\times I$}\\
\E(0.2*-0.1){(a)}\\
\E(.8*-0.1){(b)}\\
\endSetLabels 
\AffixLabels{{\includegraphics[width=4in]{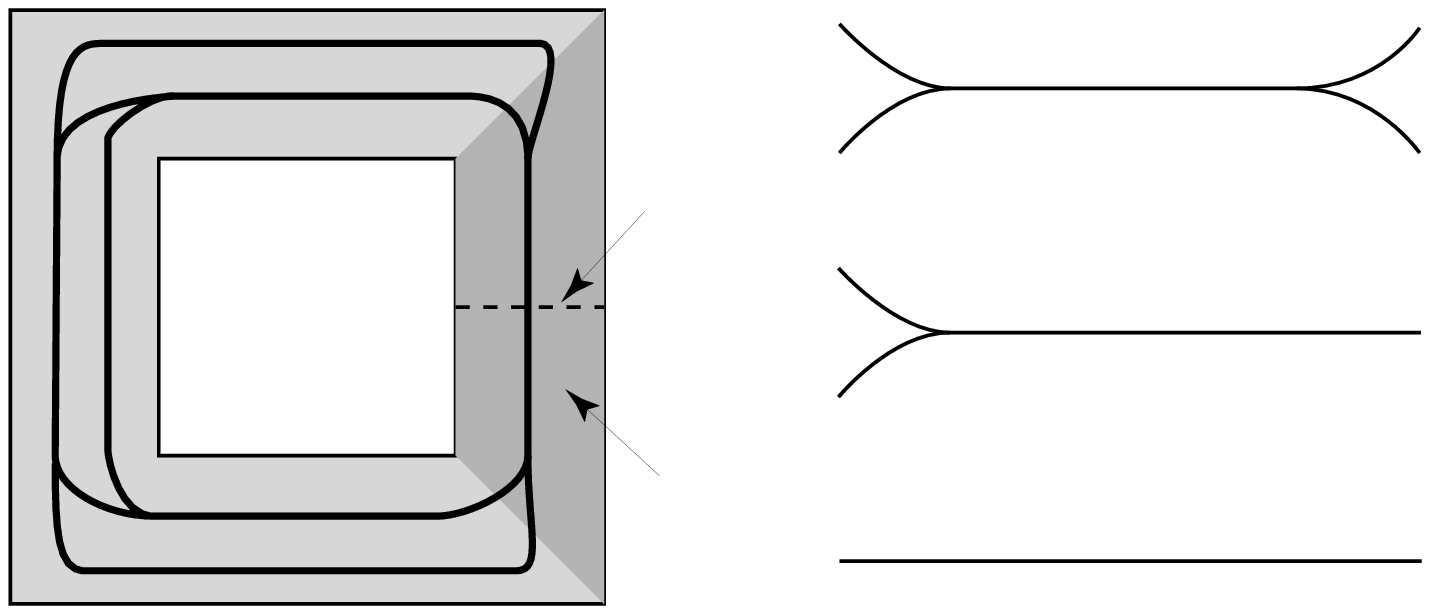}}}}
\vspace{2mm}
\caption{}\label{F33}
\end{figure}

Now, we let $D$ be a disk branch of $B$, and we identify $N_B(D)$ and
$D\times I$.  Let $O_D=\{E |\text{ $E$ is an edge of $\partial D$ with
branch direction pointing outwards.}\}$.  For each edge $E\in O_D$,
$D^B\cap (E\times I)$ (where $E\times I\subset D\times I=N_B(D)$) must
be one of the three patterns shown in Figure~\ref{F33}~(b) by our
construction.  In particular, let $p$ be the midpoint of $E$,
$\{p\}\times I\subset E\times I$ intersects $D^B\cap (E\times I)$ in a
single point.  The train track in Figure ~\ref{F33}~(a) is a picture
of the intersection of $D^B$ with $\partial D\times I$ (the shaded
annulus), where $D$ is as in Figure~\ref{F32} (a).


The next proposition is an important observation for our construction.
The notation used is the same as that in the discussion above.

\begin{proposition}\label{P31}
Let $E\in O_D$ and $p$ be the midpoint of $E$.  Then, as shown in
Figure~\ref{F33}~(a), the $I$--fiber $\{p\}\times I$ of $N_B(D)=D\times
I$ intersects $B$ in a single point.  Given any 1--dimensional
lamination $\mu$ carried by the train track $D^B\cap\partial N_B(D)$,
where $\mu\subset\partial D\times I$ is transverse to each $I$--fiber
of $\partial D\times I$, we can change the lamination $\mu$ near
$\{p\}\times I$ to get a new lamination $\mu'$ carried by
$D^B\cap\partial N_B(D)$ such that all the leaves in $\mu'$ are
circles, and hence $\mu'$ can be extended to a (2--dimensional) product
lamination carried by $D^B$.
\end{proposition}

\begin{proof}
The proof is easy.  We cut $\partial D\times I$ along $\{p\}\times
I$. Then $\mu$ is cut into a collection of compact arcs.  We can
re-glue them along $\{p\}\times I$ in such a way that these arcs close
up to become circles.  Since the intersection of $\{ p\}\times I$ and
$\partial D^B$ is a single point, the new ($1$--dimensional) lamination
is still carried by $\partial D^B$.
\end{proof}

Now we are in position to construct a lamination carried by $B$.

The first step is to construct a lamination with boundary carried by
$P(L')$.  For each branch of $P(L')$, say $S$, we construct a product
lamination $\{$\textit{a cantor set}$\}\times S$.  Since the branch
locus of $P(L')$ does not have double points, by gluing together
finitely many $\{\mathit{cantor\ sets}\}\times I$'s along the branch
locus of $P(L')$, one can easily construct a lamination fully carried
by $P(L')$.  What we want to do next is to modify this lamination so
that it can be extended to $B-int(P(L'))$, which is the union of
$D^B$'s for all the branches.

Let $D_1, D_2,\dots, D_n$ be all the disk branches of $B-L$.  Since
there is no sink disk, any disk $D_i$ has a boundary edge, say $E_i$,
with direction pointing outwards.  Locally there are $3$ branches
incident to $E_1$.  If the branch to which the branch direction of
$E_1$ points is a disk, say it is $D_2$, we denote $D_1\cup D_2$ by
$D_1\to D_2$.  Note that $E_1$ is also a boundary edge of $D_2$ with
branch direction points into $D_2$.  $D_2$ also has a boundary edge,
say $E_2$, with direction pointing outwards.  If the branch to which
$E_2$ points is a disk, say it is $D_3$, we denote $D_1\cup D_2\cup
D_3$ by $D_1\to D_2\to D_3$.  Note that the branch direction of $E_2$
points out of $D_2$ and points into $D_3$.  We proceed in this manner.
We call $D_1\to D_2\to\dots\to D_k$ a \textit{chain} if $D_i\ne D_j$
for any $i\ne j$, and call it a \textit{cycle} if $D_1=D_k$ and
$D_1\to\dots\to D_{k-1}$ is a chain.  We say that two cycles are
disjoint if there is no disk branch appearing in both cycles.  We can
decompose the union of disk branches of $B-L$ into a collection of
finitely many disjoint cycles and finitely many chains that connect
these cycles and the non-disk branches.  Moreover, we can assume that
the union of all the disk branches in those chains does not contain
any cycle, otherwise we can increase the number of disjoint cycles.
Note that a disk branch can belong to more than one chain, but it can
neither belong to more than one cycle nor belong to both a cycle and a
chain.

\begin{remark}
$D_1\to D_1$ could be a cycle.
\end{remark}

The next step is to modify the lamination and extend it through all
the chains.  For any chain $D_1\to D_2\to\dots\to D_k$, we consider
$N_B(D_i)$ and $D_i^B$, $i=1,2,\dots, k$.  There is a one-dimensional
lamination carried by $\partial D_1^B$, which is induced from the
boundary of the lamination carried by $P(L')$.  By
Proposition~\ref{P31}, we can cut $\partial D_1\times I$ along
$\{p\}\times I$ ($p\in E_1$) and re-glue it so that the
one-dimensional lamination carried $\partial D_1^B$ is a union of
circles.  There is an arc $l$ properly embedded in a branch of $P(L')$
that corresponds to $E_1$, such that one endpoint of $l$ is $p$ (in
Proposition\ref{P31}) and the other endpoint of $l$ is in $\partial
D_2^B$.  Moreover, we can cut this branch of $P(L')$ (as well as the
lamination carried by $P(L')$) along $l$, then glue it back in such a
way that the lamination carried by $P(L')$, when restricted to
$\partial D_1^B$, becomes a lamination by circles.  Clearly this
operation changes the one-dimensional lamination carried by $\partial
D_2^B$, but it does not change the boundary lamination carried by
$\partial D_i^B$ for $i\ne 1, 2$.  Now applying Proposition~\ref{P31}
to $D_2^B$, we can modify the lamination carried by $P(L')$ along
$E_2$ so that the new lamination restricted to $\partial D_2^B$ is a
lamination by circles.  Since $D_1\ne D_2$ and $D_1\to D_2\to D_1$ is
not a cycle by our assumptions, the modification along $\partial
D_2^B$ does not affect the lamination along $\partial D_1^B$, in other
words, after this modification, the lamination carried by $P(L')$
restricted to both $\partial D_1^B$ and $\partial D_2^B$ is a
lamination by circles.  We repeat this operation through the chain and
eventually get a lamination carried by $P(L')$ whose restriction to
$\partial D_i^B$ is a lamination by circles for each $i=1, 2,\dots,
k$.  We can perform this operation for all our chains and extend the
lamination of $P(L')$ through $D_i^B$ for every $D_i$ in a chain.

For any non-disk component of $B-L$, let $C$ be a simple closed curve
that is non-trivial in this component.  By Lemma~\ref{L:dl}, $C$ does
not bound a disk in $N(B)$ that is transverse to the interval fibers.
So the lamination we have constructed (for $P(L')$ and the chains) so
far does not contain any disk leaf whose boundary is in the vertical
boundary of $N_B(S)$ for any non-disk branch $S$.

By repeated application of Lemma~\ref{L26}, we can modify the
lamination and extend it through all the non-disk branches of $B$.
So, it remains to be shown that the lamination can be extended through
all the cycles.  We denote the lamination that we have constructed so
far by $\lambda$.  Note that $\lambda$ is carried by $B$ excluding
those $D_i^B$'s that correspond to the disk branches in finitely many
disjoint cycles.

\begin{figure}[ht!]
\centerline{\small
\SetLabels 
\E(0.2*-0.15){(a)}\\
\E(.8*-0.15){(b)}\\
\endSetLabels 
\AffixLabels{{\includegraphics[width=4in]{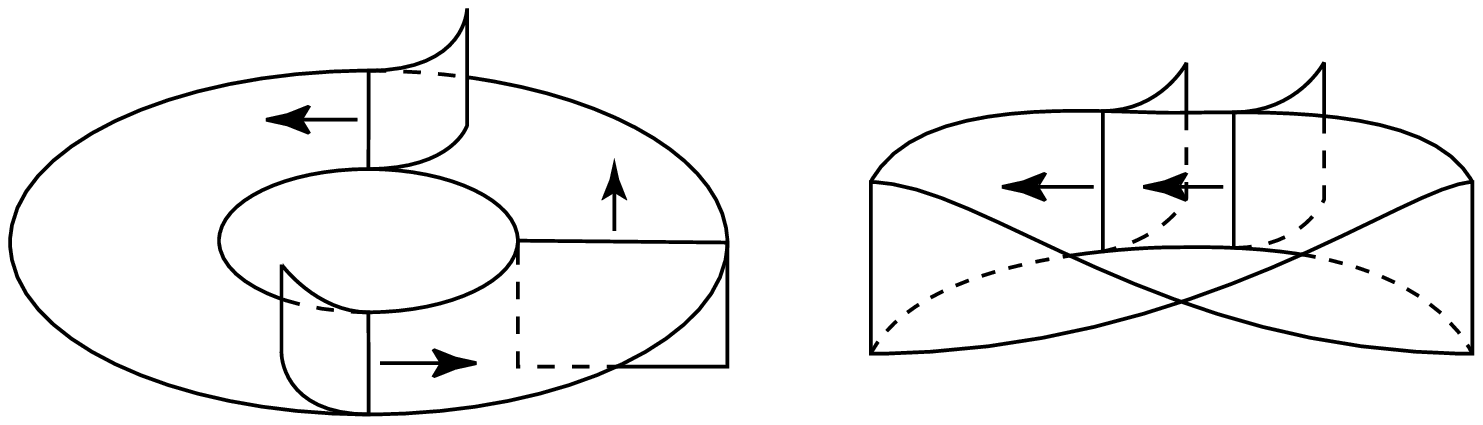}}}}
\vspace{2mm}
\caption{}\label{F34}
\end{figure}

Let $D_1\to D_2\to\dots\to D_k\to D_1$ be a cycle and $c$ be the core
of the cycle, ie, a simple closed curve in $\bigcup_{i=1}^kD_i$ such
that $c\cap D_i$ is a simple arc connecting $E_{i-1}$ to $E_i$ for
each $i$ (let $E_0=E_k$).  The intersection of $B$ with a small
regular neighborhood of $c$ in $M$, which we denote by $N(c)$, is
either a branched annulus or a branched M\"{o}bius band with coherent
branch directions, as shown in Figure~\ref{F34}.

The lamination $\lambda$ can be trivially extended to a lamination
carried by $B-\bigcup_{all\ cycles}N(c)$.  Hence, it suffices to
extend the lamination from the boundary of a branched annulus (or
M\"{o}bius band) to its interior.  We also use $\lambda$ to denote the
lamination carried by $B-\bigcup_{all\ cycles}N(c)$.

We will only discuss the branched annulus case
(ie Figure~\ref{F34}~(a)).  The branched M\"{o}bius band case is
similar.

Let $A=\bigcup_{i=1}^kd_i$ be the annulus, where $d_i=D_i\cap N(c)$,
and $A'$ be the whole branched annulus (ie, $A'=B\cap N(c)$).  Then
$\partial A=A_1\cup A_2$ consists of two circles and $A'=T_1\times
I=T_2\times I$, where $T_i$ is a train track consisting of the circle
$A_i$ and some `tails' with coherent switch directions, $i=1,2$.

Now we consider the one-dimensional lamination in $\pi^{-1}(T_i)$
induced from $\lambda$.  Since branch directions of the branched
annulus are coherent, as shown in Figure~\ref{F34} (a), the
(one-dimensional) leaves that come into $\pi^{-1}(T_i)$ from the
`tails' must be spirals with the same spiraling direction
(ie, clockwise or counterclockwise).  So the leaves coming from the
`tails' above the circle $A_i$ have the same limiting circle $H_i$,
and the leaves coming from the `tails' below the circle $A_i$ have the
same limiting circle $L_i$, $i=1,2$.  Note that the leaves may come
into $A_i$ from different sides depending on the side from which the
disks in $A'-A$ are attached to $A$, and this is what the words
`above' and `below' mean; in the branched M\"{o}bius band case, we do
not have such problems.  After replacing a leaf by an $I$--bundle over
this leaf and deleting the interior of the $I$--bundle, we can assume
that $H_i\ne L_i$, $i=1,2$.

Then we add two annuli $H$ and $L$ in $\pi^{-1}(A)$ such that
$\partial H=H_1\cup H_2$ and $\partial L=L_1\cup L_2$.  Notice that
the spirals above $H_1$ (resp. below $L_1$) in $\pi^{-1}(T_1)$ are
connected, one to one, to the spirals above $H_2$ (resp. below $L_2$)
in $\pi^{-1}(T_2)$ by the lamination $\lambda$ (restricted to
$\pi^{-1}(\partial A'-T_1-T_2)$).  So, in $\pi^{-1}(A')$, we can
naturally connect the spirals above $H_1$ (resp. below $L_1$) carried
by $T_1$ to the spirals above $H_2$ (resp. below $L_2$) carried by
$T_2$ using some ($2$--dimensional) leaves of the form $spiral\times I$
such that the boundaries of these $spiral\times I$'s lie in
$\partial\lambda$, and $H$ (resp. $L$) is the limiting annulus of
these $spiral\times I$ leaves.

There is a product region in $\pi^{-1}(A')$ between annuli $H$ and
$L$.  As in Lemma \ref{L24}, we can modify and extend the lamination
from the vertical boundary of this product region to its interior, and
hence we can extend our lamination through a given cycle.  Note that
in order to apply Lemma~\ref{L24}, we need the hypothesis that there
is no disk leaf whose boundary is in the vertical boundary of the
product region between by $H\cup L$, and this is guaranteed by
Lemma~\ref{L:ba}.

Since the cycles are disjoint by our assumption, we can successively
modify and extend the lamination through all the cycles.  The
lamination we get in the end is fully carried by $B$.  \qed

Next, we will show that if $\lambda$ is a lamination by planes, then
any branched surface that carries $\lambda$ must contain a sink disk.

Proposition~\ref{P32} was proved by Gabai (see \cite{G8}), and it is a
lamination version of a theorem of Imanishi \cite{Im} for $C^0$
foliations by planes.

\begin{proposition}\label{P32}
If M contains an essential lamination by planes, then M is
homeomorphic to the 3--torus $T^3$.
\end{proposition}

\begin{proof}
Because $\lambda$ is an essential lamination by planes, the complement
of any branched surface carrying it must be a collection of $D^2\times
I$ regions.  Hence, $\lambda$ can be trivially extended to a $C^0$
foliation by planes.  Then a theorem of Imanishi \cite{Im}, the
classification of leaf spaces and H\"{o}lder's theorem together imply
that $\pi_1(M)$ is commutative and hence $M=T^3$.
\end{proof}

\begin{proposition}\label{P33}
Let $\lambda$ be a lamination by planes in a 3--manifold $M$.  Then,
any branched surface carrying $\lambda$ cannot be a laminar branched
surface.
\end{proposition}

\begin{proof}
Suppose $B$ is a laminar branched surface that carries $\lambda$ and
$L$ is the branch locus.  If $B-L$ contains a nondisk component, then
there is a simple closed curve $C\subset B-L$ that is homotopically
nontrivial in $B-L$.  Since $B$ fully carries a lamination $\lambda$
($\lambda\subset N(B)$), $\pi^{-1}(C)\cap\lambda$ must contain a
simple closed curved in a leaf of $\lambda$.  Since every leaf of
$\lambda$ is a plane, this simple closed curve bounds an embedded disk
$D$ in a leaf.  The disk clearly satisfies the two conditions in
Lemma~\ref{L:dl} with $\pi(\partial D)=C$, which contradicts the
assumption that $B$ is a laminar branched surface.  Therefore, $B-L$
is a union of disks.  Since there is no sink disk, every disk branch
of $B$ has an edge whose branch direction points outwards.  Then the
disk branches of $B$ form at least one cycle as before. A small
neighborhood of the core of the cycle is either a branched annulus or
a branched M\"{o}bius band, as show in Figure~\ref{F34}.  Hence, one
gets a curve with non-trivial holonomy, which contradicts the Reeb
Stability Theorem and the assumption that every leaf is a plane.
\end{proof}

\section{Splitting branched surfaces along laminations}

Suppose $\lambda$ is an essential lamination in an orientable
3--manifold $M$ and $\lambda$ is not a lamination by planes.  By
`blowing air' into leaves, ie, replacing every leaf by an $I$--bundle
over this leaf and then deleting the interior of the $I$--bundle, we
can assume that $\lambda$ is nowhere dense and fully carried by a
branched surface $B$.  By \cite{GO}, we can assume that $B$ satisfies
the conditions in Proposition~\ref{P11}.  It is easy to see that $B$
still satisfies conditions 1, 2, and 3 in Proposition~\ref{P11} after
any further splitting along $\lambda$.  We will show in this section
that we can split $B$ along $\lambda$ to make it a laminar branched
surface.

Let $B'$ be a union of some branches of $B$.  We call a point of $B'$
an interior point if it has a small open neighborhood in $M$ whose
intersection with $B'$ is a small branched surface without boundary,
ie, the intersection is one of the 3 pictures as shown in
Figure~\ref{F11}, otherwise we call it a boundary point.  We denote
the union of boundary points of $B'$ by $\partial B'$.  Next, we give
every arc in $\partial B'$ a normal direction pointing into $B'$, and
give every arc in $L\cap (B'-\partial B')$ its branch direction.  We
call the direction that we just defined for $\partial B'$ and $L\cap
(B'-\partial B')$ the \textit{direction associated with $B'$}.  We
call $B'$ (and also $N(B')=\pi^{-1}(B')$) a \textit{safe region} if it
satisfies the following conditions:

\begin{enumerate}
\item $B'$ does not contain any disk branch with the induced direction
(from the direction associated with $B'$ that we just defined) of
every boundary arc pointing inwards;
\item for any non-disk branch in $B'$, if the direction (associated
with $B'$) of every boundary arc points inwards, then it contains a
closed curve that is homotopically non-trivial in $M$.
\end{enumerate}

Thus, by our definition, every disk branch (of $B$) lying in a safe
region $B'$ must have a boundary arc lying in the interior of $B'$
with branch direction pointing outwards.

\begin{proposition}\label{P41}
Let $B'$ be a safe region.  For any smooth arc $\alpha\subset L$, if
either $\alpha\subset B'-\partial B'$ or $\alpha\subset\partial B'$
and the branch direction of $\alpha$ points into $B'$, then the union
of $B'$ and all the branches that are incident to $\alpha$ is still a
safe region.
\end{proposition}

\begin{proof}
Let $D\nsubseteq B'$ be a branch incident to $\alpha$.  Then the
branch direction of $\alpha$ points out of $D$.  Hence $B'\cup D$
still satisfies our conditions for safe regions.
\end{proof}

\begin{proposition}\label{P42}
Let $B'$ be a safe region.  If $B'=B$, then the branched surface $B$
is a laminar branched surface.\qed
\end{proposition}

\begin{figure}[ht!]
\centerline{\small
\SetLabels 
\E(.5*0.57){(1)}\\
\E(.49*0.29){(2)}\\
\E(.5*0.41){splitting}\\
\endSetLabels 
\AffixLabels{{\includegraphics[width=4in]{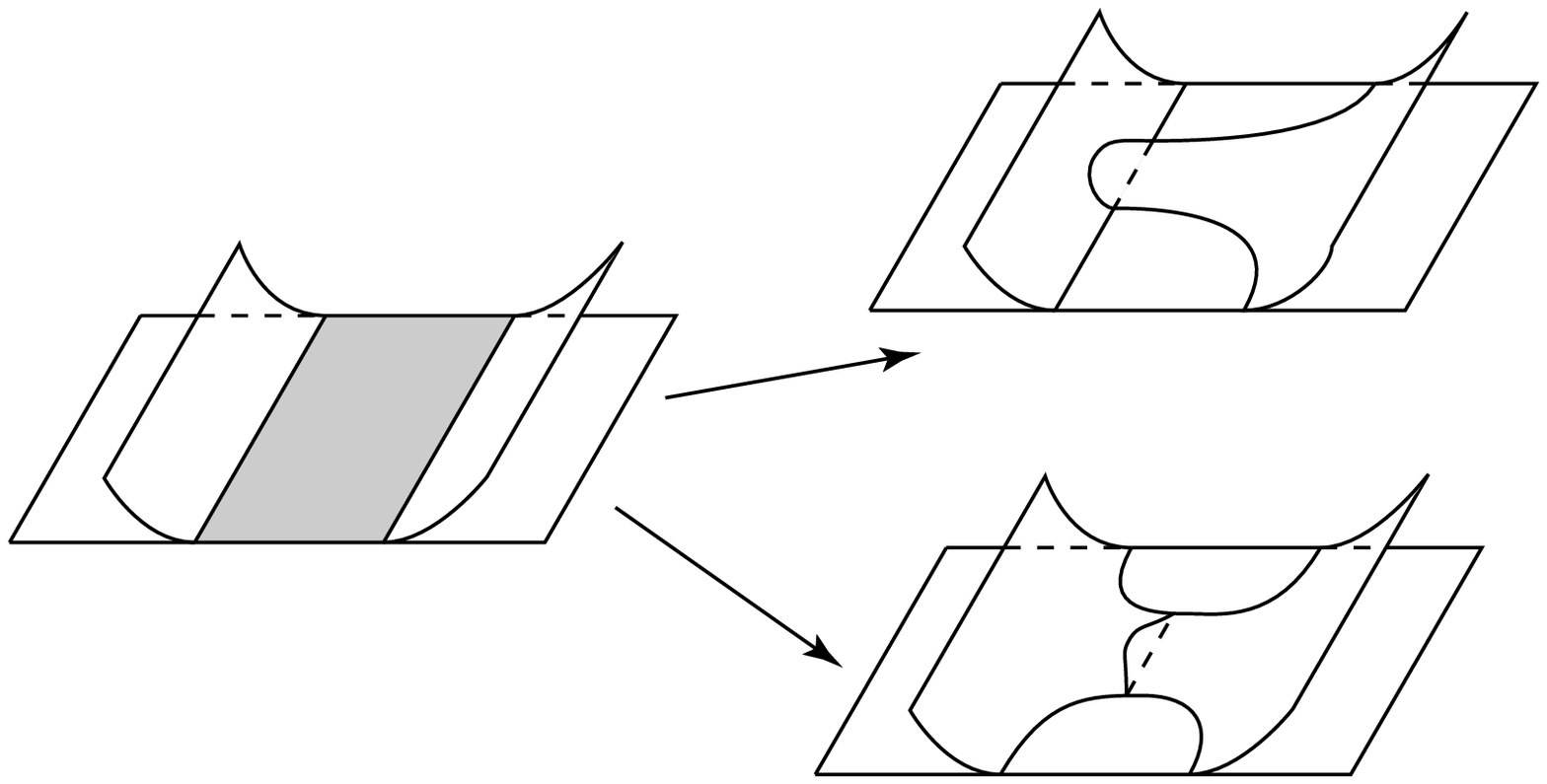}}}}
\vspace{2mm}
\caption{}\label{F41}
\end{figure}

\begin{figure}[ht!]
\centerline{\small
\SetLabels 
\E(.125*0.67){$B'$}\\
\E(.435*0.67){$S$}\\
\E(.68*0.37){$B'$}\\
\E(.97*1){$S$}\\
\E(.4*0.815){safe region}\\
\E(.52*0.545){(1)}\\
\E(.52*0.29){(2)}\\
\E(.5*0.395){splitting}\\
\endSetLabels 
\AffixLabels{{\includegraphics[width=4in]{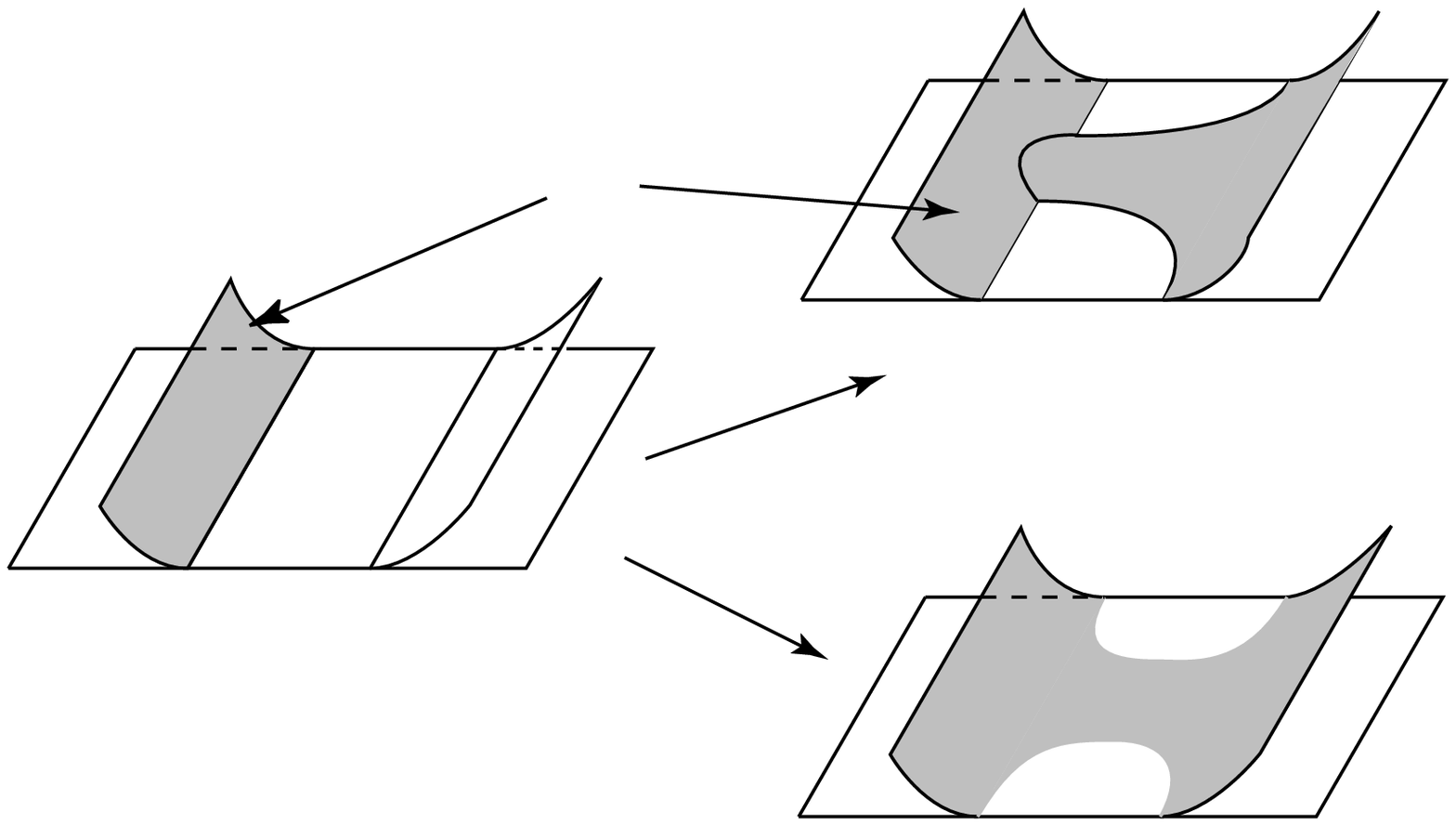}}}}
\vspace{2mm}
\caption{}\label{F42}
\end{figure}

Next, we will show how the safe region changes when we split the
branched surface.  What we want is to enlarge our safe region by
splitting the branched surface.  Suppose we do some splitting to $B$
whose local picture is shown in Figure~\ref{F41}.  We will call the
splitting an \textit{unnecessary splitting} if the shaded area in
Figure~\ref{F41} belongs to the safe region, otherwise, we call it a
\textit{necessary splitting}.  Note that, by Proposition~\ref{P41}, if
the shaded area belongs to the safe region, we can include all the
branches (of the branched surface on the left) in Figure~\ref{F41}
into the safe region, so we do not need to do such a splitting to
enlarge our safe region.  The following proposition says that the safe
region does not decrease under necessary splittings.  The proof is an
easy application of Proposition~\ref{P41}.  Figure~\ref{F42} is just
Figure~\ref{F41} with different shaded regions which denote the safe
region.

\begin{proposition}\label{P43}
Let $B'$ be a safe region.  If after a necessary splitting, a branch
$S$ of $B$ slides on $B'$, as shown in splitting (1) in
Figure~\ref{F42}, or $S$ and a branch in $B'$ locally becomes one
branch, as shown in splitting (2) in Figure~\ref{F42}, then we can
enlarge the safe region after the splitting as shown in the two
pictures on the right in Figure~\ref{F42}.  In particular, for any
interval fiber of $N(B)$ that is in a safe region, if this fiber
breaks into two interval fibers after some necessary splitting, then
we can enlarge the safe region after the splitting such that both
interval fibers lie in the safe region.
\end{proposition}

\begin{proof}
After the splitting (1) in Figure~\ref{F42}, $S$ has a boundary arc,
with branch direction pointing outwards, lying in the interior of the
shaded region in Figure~\ref{F42}.  So, by Proposition \ref{P41},
after splitting (1), the union of the original safe region and this
branch is still a safe region.  Note that, since it is a necessary
splitting, the branch in the middle does not belong to the safe region
and the change (done by the splitting) of this branch does not affect
the safe region.

In the splitting (2) of Figure~\ref{F42}, if $S$ and the shaded region
in the left picture of Figure~\ref{F42} do not belong to the same
branch in $B$, then after splitting (2) the new shaded branch in $B$
either has a boundary arc with direction (associated with the safe
region) pointing outwards or contains a non-trivial curve, since the
shaded branch before the splitting is in the safe region.

Suppose $S$ and the shaded region in the left picture of
Figure~\ref{F42} belong to the same branch in $B$ (before the
splitting).  We denote this branch by $D$ ($D\subset B'$).  If $D$ has
a boundary arc whose direction (associated with the safe region)
points outwards, then after the splitting, it still has such a
boundary arc.  If $D$ does not have such an arc, then $D$ is not a
disk and it contains a curve that is homotopically non-trivial in $M$.
Thus, after splitting (2), the branch still contains such an essential
closed curve, and we can include this branch (after the splitting)
into the safe region.
\end{proof}

Now we are ready to prove the following lemma that is a half of
Theorem~\ref{T01}.  In the proof, we first construct a safe region
$B'$ by taking a small neighborhood of a union of finitely many
essential curves in leaves of $\lambda$.  Then, we perform some
necessary splitting (along $\lambda$) and enlarge $B'$ so that $B-B'$
lies in a union of disjoint 3--balls.  After splitting $B$ along the
boundary of these 3--balls and getting rid of the disks of contact in
these 3--balls, we can include the whole of $B$ to the safe region, and
hence $B$ becomes a laminar branched surface.

\begin {lemma}\label{L44}
Let $\lambda$ be an essential lamination that is not a lamination by
planes.  Then $\lambda$ is carried by a laminar branched surface.
\end{lemma}

\begin{proof}
Suppose that $\lambda$ is an essential lamination in a 3--manifold $M$.
We first show that any sub-lamination $\mu$ of $\lambda$ is not a
lamination by planes.  Suppose $\mu$ is a lamination by planes and
$\gamma$ is a closed curve in $M-\mu$.  Then, by splitting $\mu$, we
can assume that $\mu$ is carried by $N(B)$ and $\gamma$ lies in a
component $C$ of $M-int(N(B))$.  By isotoping $\mu$, we can assume
that $\partial_hN(B)\subset\mu$.  Let $l$ be a boundary leaf of the
component of $M-\mu$ that contains $\gamma$.  Then, we can choose a
big disk $D_1\subset l$ such that $l\cap\partial_hN(B)\subset D_1$.
Moreover, there is a vertical annulus $A$ consisting of subarcs of
$I$--fibers of $N(B)$ such that one boundary circle of the $A$ is
$\partial D_1$, $\partial A\subset\mu$, and the $int(A)$ lies in the
same component of $M-\mu$ that contains $\gamma$.  Since $\mu$ is
assumed to be a lamination by planes, the other boundary circle of the
$A$ bounds a disk $D_2$ in a leaf of $\mu$.  So, $D_1\cup A\cup D_2$
forms a sphere.  Since $M$ is irreducible, $D_1\cup A\cup D_2$ bounds
a 3--ball whose interior lies in $M-\mu$.  As
$l\cap\partial_hN(B)\subset D_1$, $C$ and hence $\gamma$ must lie in
this 3--ball, which implies every component of $M-\mu$ is simply
connected.  Since every leaf of $\lambda$ is $\pi_1$--injective, every
leaf of $\lambda$ must be a plane, which contradicts our hypothesis.
Therefore, any sub-lamination of $\lambda$ is not a lamination by
planes.

By \cite{GO}, we can assume that $\lambda$ is carried by a branched
surface $B$ that satisfies the conditions in Proposition~\ref{P11}.
Moreover, we also assume that $\lambda$ is in Kneser-Haken normal form
with respect to a triangulation $\mathcal{T}$ (one can take a fine
enough triangulation so that the branched surface $B$ is a union of
normal disks in this triangulation).  Now, $\lambda$ lies in $N(B)$
transversely intersecting every interval fiber of $N(B)$, and
$N(B)\cap\mathcal{T}^{(1)}$ is a union of $I$--fibers of $N(B)$, where
$\mathcal{T}^{(1)}$ is the $1$--skeleton of $\mathcal{T}$.  Note that,
by \cite{GO}, $B$ still satisfies conditions 1--3 in
Proposition~\ref{P11} after any splitting.

For any point $x\in\lambda\cap\mathcal{T}^{(1)}$, we denote the leaf
that contains $x$ by $l_x$.  Let $\mu_x$ be the closure of $l_x$.
Then $\mu_x$ is a sub-lamination of $\lambda$.  Hence, $\mu_x$ is not
a lamination by planes.  Let $c_x$ be a non-trivial simple closed
curve in a non-plane leaf of $\mu_x$.  Then there is an embedding $A\co 
S^1\times I\to N(B)$, where $I=[-1,1]$, such that $A(\{ p\}\times I)$
is a sub-arc of an interval fiber of $N(B)$, $A(S^1\times\{ 0\}
)=c_x$, and every closed curve in $A^{-1}(\lambda )$ is of the form
$S^1\times\{ t\}$ for some $t\in (-1,1)$.  Moreover, after some
isotopies, we can assume that there are $a,b\in I$ such that $-1<a\le
0\le b<1$, $A(S^1\times\{ a,b\})\subset\lambda$, and $A^{-1}(\lambda
)\cap (S^1\times (I-[a,b]))$ is either empty or a union of spirals
whose limiting circles are $S^1\times\{ a,b\}$.  We call such embedded
annuli \textit{regular annuli}.

To simplify notation, we will not distinguish the map $A$ and its
image.  Since $c_x$ lies in the closure of $l_x$, there must be a
simple arc in $l_x$ connecting $x$ to the annulus $A$.  Moreover,
there is an embedding $b\co  I\times (-1,1)\to N(B)$ such that
$b(-1,0)=x$, $b(\{ 1\}\times (-1,1))\subset A$, $I_x=b(\{ -1\}\times
(-1,1))\subset\mathcal{T}^{(1)}$, and $b^{-1}(\lambda )$ is a union of
compact parallel arcs connecting $\{ -1\}\times (-1,1)$ to $\{
1\}\times (-1,1)$.  Hence, $I_x$ is an open neighborhood of $x$ in
$\mathcal{T}^{(1)}$.  For every point
$x\in\lambda\cap\mathcal{T}^{(1)}$, we have such an open interval
$I_x$ and a regular annulus as $A$ above.  By compactness, there are
finitely many points $x_1, x_2,\dots, x_n$ in
$\lambda\cap\mathcal{T}^{(1)}$ such that
$\lambda\cap\mathcal{T}^{(1)}\subset\cup_{i=1}^nI_{x_i}$, where
$I_{x_i}$ is an open neighborhood of $x_i$ in $\mathcal{T}^{(1)}$ as
above.  Let $A_1,\dots, A_n$ be the regular annuli that correspond to
$x_1, \dots, x_n$ respectively as above.  Since
$\lambda\cap\mathcal{T}^{(1)}\subset\cup_{i=1}^nI_{x_i}$, for any
$x\in\lambda\cap\mathcal{T}^{(1)}$, there is an arc on a leaf of
$\lambda$ connecting $x$ to $A_i$ for some $i$.  Note that each $A_i$
is embedded but $A_i$ and $A_j$ may intersect each other if $i\ne j$.

\medskip
\noindent
\textbf{Claim}\qua  There are finitely many disjoint regular annuli
$E_1,\dots, E_k$ such that, for any
$x\in\lambda\cap\mathcal{T}^{(1)}$, there is an arc in a leaf of
$\lambda$ connecting $x$ to $E_i$ for some $i$.

\begin{proof}[Proof of the Claim]
If $A_1,\dots, A_n$ are disjoint, the claim holds immediately.
Suppose that $A_1\cap A_2\ne\emptyset$.  As a map, $A_i\co S^1\times I\to
N(B)$ is an embedding ($i=1,\dots, n$).  To simplify notation, we use
$A_i$ to denote both the map and its image in $N(B)$.  After some
homotopies, we can assume that $A_i^{-1}(A_j)$ is a union of disjoint
sub-arcs of the $I$--fibers of $S^1\times I$, and $(S^1\times\partial
I)\cap A_i^{-1}(\lambda )\cap A_i^{-1}(A_j)=\emptyset$ ($i\ne j$).
Thus, the intersection of $A_i^{-1}(\lambda )$ and $A_i^{-1}(A_j)$
must lie in the interior of $A_i^{-1}(A_j)$.

$A_1^{-1}(\lambda )$ is a one-dimensional lamination in $S^1\times I$,
and by our construction, every leaf of $A_1^{-1}(\lambda )$ that is
not a circle must have a limiting circle in $A_1^{-1}(\lambda )$.
Therefore, if every circular leaf of $A_1^{-1}(\lambda )$ has
non-empty intersection with $A_1^{-1}(\cup_{i=2}^nA_i)$, then for
every point $p\in A_1(A_1^{-1}(\lambda ))$, there is an arc in a leaf
of $\lambda$ connecting $p$ to $\cup_{i=2}^nA_i$.  Hence, for any
point $x\in\lambda\cap\mathcal{T}^{(1)}$, there is an arc in $l_x$
connecting $x$ to $\cup_{i=2}^nA_i$, and we only need to consider
$n-1$ annuli $A_2, \dots, A_n$.  If there are circular leaves in
$A_1^{-1}(\lambda )$ whose intersection with
$A_1^{-1}(\cup_{i=2}^nA_i)$ is empty, then since $A_1^{-1}(\lambda
)\cap A_1^{-1}(\cup_{i=2}^nA_i)$ lies in the interior of
$A_1^{-1}(\cup_{i=2}^nA_i)$, there are finitely many disjoint annuli
$B_1,\dots, B_m$ in $S^1\times I$ such that
$A_1^{-1}(\cup_{i=2}^nA_i)\cap (\cup_{j=1}^mB_j)=\emptyset$, every
circular leaf of $A_1^{-1}(\lambda )$ either has non-empty
intersection with $A_1^{-1}(\cup_{i=2}^nA_i)$ or lies in $B_j$ for
some $j$, and $A_1|_{B_j}$ is a regular annulus for each $j$.  To
simplify notation, we will not distinguish $A_1|_{B_j}$ and its image
in $N(B)$.  Thus, for any point $x\in\lambda\cap\mathcal{T}^{(1)}$,
there is an arc in a leaf of $\lambda$ connecting $x$ to
$(\cup_{i=2}^nA_i)\cup (\cup_{j=1}^mA_1|_{B_j})$, and $A_1|_{B_j}$ is
disjoint from $A_i$ for any $i,j$ ($i\ne 1$).

By repeating the construction above, eventually we will get finitely
many such disjoint regular annuli as in the claim.
\end{proof}

As in the claim, $E_1,\dots, E_k$ are disjoint regular annuli.  Let
$N(E_i)$ be a small neighborhood of $E_i$ in $M$ such that $N(E_i)\cap
N(E_j)=\emptyset$ if $i\ne j$.  Topologically, $N(E_i)$ is a solid
torus for each $i$.  $N(E_i)\cap\lambda$ consists of a union of
parallel annuli and some simply connected leaves.  The limit of every
simply connected leaf in $N(E_i)\cap\lambda$ is either an annulus or a
union of two annuli depending on the number of ends of the leaf.
Since the $N(E_i)$'s are disjoint, by `blowing air' into the leaves,
we can split the branched surface $B$ along $\lambda$ such that every
component of $B\cap N(E_i)$ is either an annulus or a branched annulus
with coherent branch directions, as shown in Figure~\ref{F34}~(a),
whose core is homotopically essential in $M$.  Let $D$ be a component
of $B\cap N(E_i)$ that is a branched annulus.  Since $D$ has coherent
branch directions, every branch in $D$ has a boundary edge with branch
direction pointing outwards.  Since the core of every solid torus
$N(E_i)$ is an essential curve in $M$, the union of the branches of
$B$ that have non-empty intersection with $\cup_{i=1}^kN(E_i)$ is a
safe region.  We denote the safe region by $B'$ and
$N(B')=\pi^{-1}(B')$ as before.

Note that if $B$ contains a trivial bubble, then we can collapse the
trivial bubble without destroying the branched annuli constructed
above, though the number of ``tails" in a branched annulus may
decrease.  More precisely, let $c$ be the core of a branched annulus
as above, by the definition of trivial bubble, we can always pinch $B$
to eliminate a trivial bubble so that the neighborhood of $c$ after
this pinching is still a branched annulus with coherent branch
direction.  Moreover, if $B\cap N(E_i)$ is an annulus, since the core
of $N(E_i)$ is homotopically nontrivial, the operation of eliminating
trivial bubbles does not affect the annulus $B\cap N(E_i)$.  Thus, our
safe region will never be empty due to eliminating trivial bubbles.
Next, we will perform necessary splitting to our branched surface.  If
we see any trivial bubble during the splitting, we eliminate it by
pinching the branched surface and start over.  Since the number of
components of the complement of the branched surface never increases
during necessary splittings and the number of components decreases by
one after a trivial bubble is eliminated, eventually we will never get
any trivial bubble.  Therefore, we can assume the necessary splittings
we perform in the following never create any trivial bubble.

For any $x\in\lambda\cap\mathcal{T}^{(1)}$, by the claim and our
construction of $B'$, there is a simple arc $\gamma\co  [0,1]\to l_x$
connecting $x$ to a point in $N(B')$, ie, $\gamma (0)=x$ and $\gamma
(1)=y\in N(B')$.  Moreover, there is an embedding $b_x\co  [0,1]\times
(-\epsilon, \epsilon )\to N(B)$ such that
$b_x|_{[0,1]\times\{0\}}=\gamma$, $b_x(\{t\}\times (-\epsilon,
\epsilon ))$ is a subarc of an $I$--fiber of $N(B)$, and
$b_x^{-1}(\lambda)$ is a union of parallel arcs connecting
$\{0\}\times (-\epsilon, \epsilon)$ to $\{1\}\times (-\epsilon,
\epsilon)$.  For each $t\in [0,1]$, we denote $b_x(\{t\}\times
(-\epsilon, \epsilon ))$ by $I_t$, and we also use $I_x$ to denote
$I_0$.  Thus, $x\in I_x$.  To simplify notation, we do not distinguish
$b_x$ and its image in $N(B)$.  So, $b_x=\cup_{t\in [0,1]}I_t$.  Now
for every $x\in\lambda\cap\mathcal{T}^{(1)}$, we have such an open
interval $I_x\subset N(B)\cap\mathcal{T}^{(1)}$ as above.  By
compactness, we can choose finitely many points $x_1, x_2,\dots,
x_n\in\lambda\cap\mathcal{T}^{(1)}$ such that
$\lambda\cap\mathcal{T}^{(1)}\subset\cup_{i=1}^nI_{x_i}$.

Using Proposition~\ref{P41}, we enlarge our safe region as much as we
can.  Then we do the necessary splitting along $\lambda$, and include
all possible branches into our safe region (using
Proposition~\ref{P43}) after the splitting.  We will always denote the
safe region by $B'$ (or $N(B')$).  If a certain splitting cuts through
a band $b_x$ and a vertical arc $I_t$ of $b_x$ breaks into some
smaller arcs $J_{t_1}, J_{t_2},\dots, J_{t_h}$, to simplify notation,
we will denote $\cup_{t=1}^hJ_{t_i}$ also by $I_t$, and denote
$\cup_{t\in [0,1]}I_t$ also by $b_x$.  Note that our new bands after
splitting may contain some `bubbles', as shown in Figure~\ref{F43},
but they are always embedded in $N(B)$ by our construction.

\begin{figure}[ht!]
\vspace{2mm}
\centerline{\small
\SetLabels 
\E(.43*1.08){band $b_x$}\\
\E(.5*0.4){splitting}\\
\endSetLabels 
\AffixLabels{{\includegraphics[width=4in]{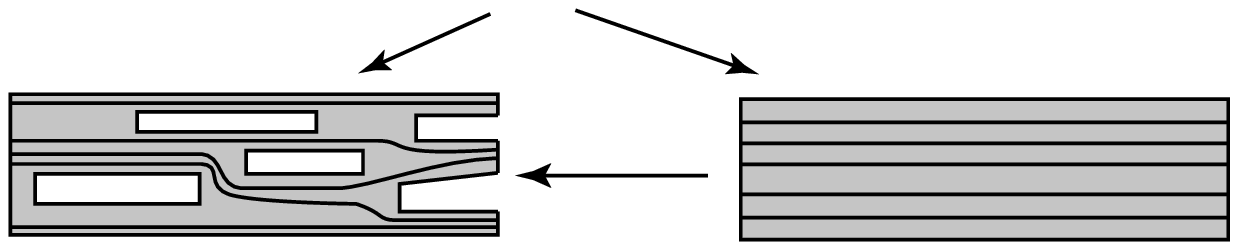}}}}
\vspace{2mm}
\caption{}\label{F43}
\end{figure}

Next, we will show that we can do some necessary splitting along a
band $b_x$ and include $I_x$ in the safe region.  By
Proposition~\ref{P43}, we know that once an interval fiber of $N(B)$
is in the safe region, it will stay in the safe region forever, though
it may break into some small intervals after further splitting.

Suppose after some splitting, the band $b_x$ contains some `bubbles'
as shown in Figure~\ref{F43}.  Although $b_x$ is embedded in $N(B)$,
there may be an $I$--fiber of $N(B)$ whose intersection with $b_x$ has
more than one component.  By perturbing $b_x$ a little, we can assume
that there are only finitely many $I$--fibers of $N(B)$ whose
intersection with $b_x$ have more than one component.  So, $\pi (b_x)$
is an immersed train track (immersed curve with some `bubbles') on
$B$, where $\pi$ is the map collapsing every interval fiber to a
point.  Those finitely many $I$--fibers whose intersection with $b_x$
have more than one component become double points of $\pi(b_x)$ after
the collapsing.  Let $\mathcal{C}$ be the number of double points of
$\pi (b_x)-B'$.  If $\mathcal{C}=0$, then we perform all possible
necessary splittings along $b_x$ and enlarge our safe region as in
Proposition \ref{P43}.  Since $b_x$ is compact, after finitely many
necessary splittings along $b_x$, the whole of $b_x$ and (hence $I_x$)
is included in the safe region.

Let $\cup_{t\in (a,b)}I_t\subset b_x-N(B')$ and $I_a\subset N(B')$
($[a,b]\subset [0,1]$).  Suppose that $\pi (\cup_{t\in (a,b)}I_t)$
contains double points.  We split $N(B)$ between $I_a$ and $I_b$ along
$b_x$ (using only necessary splittings) as above.  After the splitting
passes an interval fiber that is the inverse image (ie, $\pi^{-1}$)
of a double point of $\pi (b_x)-B'$, either the double point
disappears under the collapsing map of the new branched surface after
the splitting, or it is included in the safe region, as shown in
Figure~\ref{F44}.  Therefore, $\mathcal{C}$ decreases and eventually
we can include the whole band $b_x$ in the safe region.

\begin{figure}[ht!]
\vspace{2mm}
\centerline{\small
\SetLabels 
\E(.3*.75){band $b_x$}\\
\E(.57*0.5){splitting}\\
\E(.4*0.9){safe region}\\
\E(.215*0.2){safe region}\\
\endSetLabels 
\AffixLabels{{\includegraphics[width=4in]{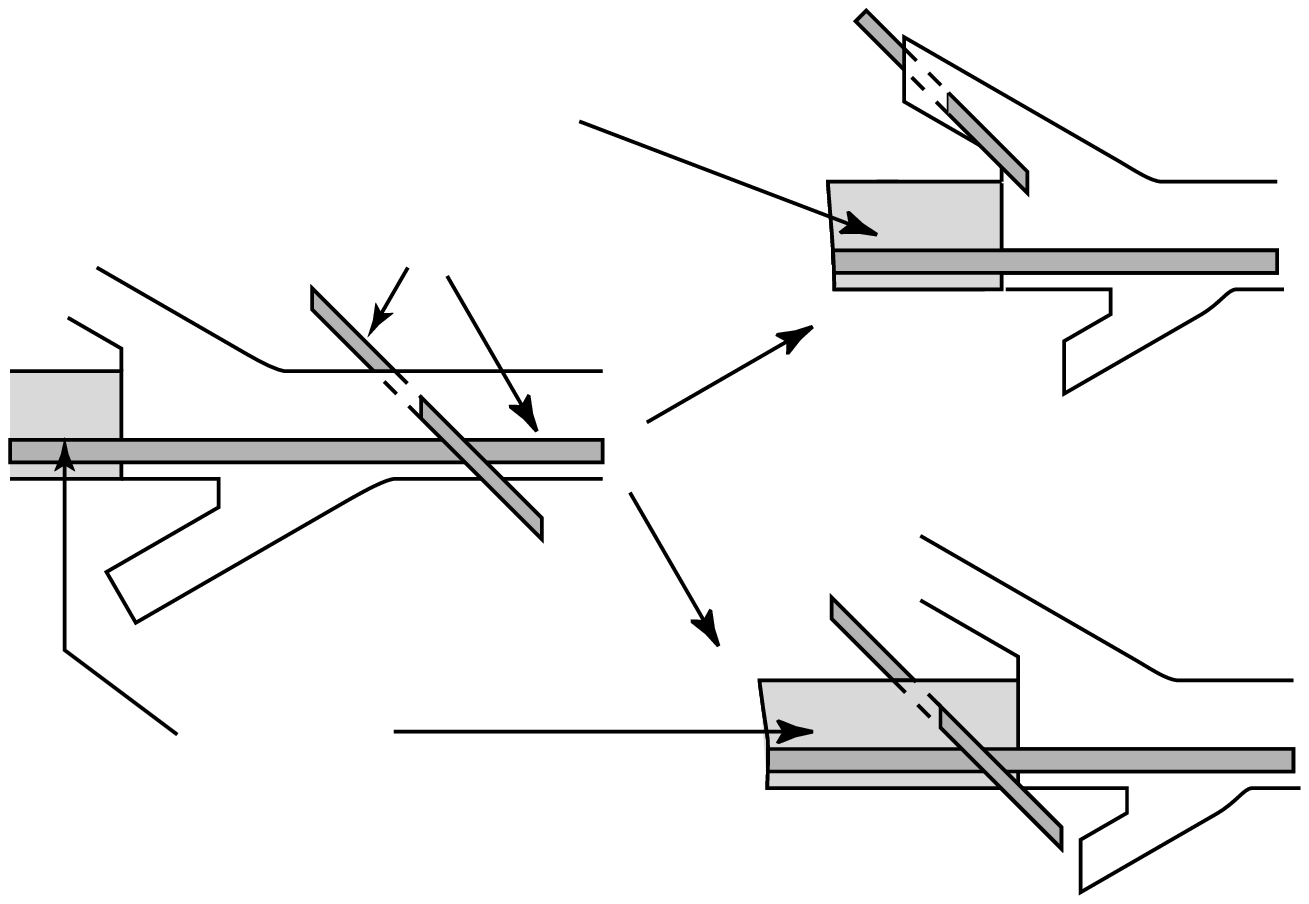}}}}
\vspace{2mm}
\caption{}\label{F44}
\end{figure}

Since there are finitely many such intervals that cover
$\lambda\cap\mathcal{T}^{(1)}$, we can include
$N(B)\cap\mathcal{T}^{(1)}$ into the safe region after finitely many
steps.  Then, by performing similar splittings, we can include
$N(B)\cap\mathcal{T}^{(2)}$ into the safe region.  Now $N(B)-N(B')$ is
contained in the interior of finitely many disjoint $3$--simplices,
ie $3$--balls.

We consider $B-(B'-\partial B')$, and let $\Gamma_1, \dots,\Gamma_s$
be the components of $B-(B'-\partial B')$.  Each $\Gamma_i$ is a union
of branches of $B$.  We can define the boundary of $\Gamma_i$ in the
same way as we did for $B'$ at the beginning of this section.  The
branch direction of every boundary arc of any $\Gamma_i$ must point
into $\Gamma_i$, and the other two (local) branches incident to this
arc must belong to $B'$ (since it is a boundary arc of $\Gamma_i$),
otherwise, using Proposition~\ref{P41}, we can enlarge $B'$ by adding
all the branches incident to this arc to $B'$.  Thus, for each
$\Gamma_i$, there is a small neighborhood of $\Gamma_i$, which we
denote by $N(\Gamma_i)$, such that $N(\Gamma_i)\cap
N(\Gamma_j)=\emptyset$ if $i\ne j$.  Moreover, after some necessary
splitting, we can assume that each $N(\Gamma_i)$ is homeomorphic to a
$3$--ball.  By our definition of the safe region, any branch in $B'$
that has non-empty intersection with $\cup_{i=1}^s\partial
N(\Gamma_i)$ either contains an essential closed curve, or has a
boundary arc (lying in the interior of $B'$) with branch direction
pointing outwards.  Therefore, after any (unnecessary) splitting along
$B\cap\partial N(\Gamma_i)$, each branch in $B-int(N(\Gamma_i))$
either contains an essential closed curve, or has a boundary arc
(lying in the interior of $B-int(N(\Gamma_i))$) with branch direction
pointing outwards.  Next, we split $B$ along $\lambda\cap\partial
N(\Gamma_i)$ (for each $i$) so that $B\cap (\cup_{i=1}^s\partial
N(\Gamma_i))$ becomes a union of circles, and at this point, each
branch of $B-int(N(\Gamma_i))$ that has non-empty intersection with
$\partial N(\Gamma_i)$ either contains an essential closed curve, or
has a boundary arc (lying in the interior of $B-int(N(\Gamma_i))$)
with branch direction pointing outwards.  Then, we split $B$ to get
rid of the disks of contact in $\cup_{i=1}^s int(N(\Gamma_i))$.  After
this splitting, $B\cap N(\Gamma_i)$ becomes a union of disks for each
$i$, and each branch of $B$ either contains an essential closed curve,
or has a boundary arc with branch direction pointing outwards.  Hence,
$B$ contains no sink disk after all these splittings, and becomes a
laminar branched surface.
\end{proof}


\begin{thebibliography}
\bibitem{D} \textbf{C Delman}, \emph{Essential laminations and Dehn
surgery on $2$--bridge knots}, Topology Appl. 63 (1995)
201--221

\bibitem{EHN} \textbf{D Eisenbud{\rm,} U Hirsch{\rm,} W Neumann},
\emph{Transverse foliations of Seifert bundles and self-homeomorphism
of the circle}, Comment. Math. Helv. 56 (1981) 638--660

\bibitem{G1} \textbf{D Gabai}, \emph{Foliations and the topology of
3--manifolds}, J. Differential Geometry, 18 (1983) 445--503

\bibitem{G2} \textbf{D Gabai},  \emph{Foliations and the topology of
3--manifolds II}, J. Differential Geometry, 26 (1987) 461--478

\bibitem{G3} \textbf{D Gabai}, \emph{Foliations and the topology of
3--manifolds III}, J. Differential Geometry, 26 (1987) 479--536

\bibitem{G4} \textbf{D Gabai}, \emph{Taut foliations of 3--manifolds
and suspensions of $S^1$}, Ann. Inst. Fourier, Grenoble, 42 (1992)
193--208

\bibitem{G5} \textbf{D Gabai},  \emph{Problems in Foliations and
Laminations}, Stud. in Adv. Math. AMS/IP, 2 (1997) 1--34

\bibitem{G8} \textbf{D Gabai},  \emph{Foliations and
3--manifolds}, Proceedings of the International Con\-gress of
Mathematicians, Vol. I (Kyoto, 1990) 609--619

\bibitem{GO} \textbf{D Gabai}, \textbf{U Oertel}, \emph{Essential
laminations in 3--manifolds}, Ann. of Math. 2 (1989) 41--73

\bibitem{GK} \textbf{D Gabai}, \textbf{W\,H Kazez}, \emph{Group
negative curvature for 3--manifolds with genuine laminations},
Geometry and Topology, 2 (1998) 65--77

\bibitem{HRS} \textbf{J Hass{\rm,} H Rubinstein{\rm,} P Scott},
\emph{Compactifying coverings of closed 3--manifolds}, J. Differential
Geomerry,  30 (1989) 817--832

\bibitem{Ha1} \textbf{A Hatcher},  \emph{On the Boundary Curves of
Incompressible Surfaces}, Pacific J. Math.  99 (1982) 373--377

\bibitem{Ha2} \textbf{A Hatcher},  \emph{Some examples of essential
laminations in 3--manifolds}, Ann. Inst. Fourier, Grenoble, 42 (1992)
313--325

\bibitem{Im} \textbf{H Imanishi},  \emph{On the theorem of
Denjoy-Sacksteder for codimension one foliations without holonomy},
J. Math. Kyoto Univ. 14 (1974) 607--634

\bibitem{L} \textbf{T Li}, \emph{Laminar branched surface in
3--manifolds II}, in preparation

\bibitem{Na} \textbf{R Naimi},  \emph{Constructing essential
laminations in 2--bridge knot surgered 3--manifolds},  Pacific J. Math.
180 (1997) 153--186

\bibitem{N} \textbf{S Novikov},  \emph{Topology of Foliations},
Moscow Math. Soc.  14 (1963) 268--305

\bibitem{PH} \textbf{R\,C Penner{\rm,} J\,L Harer},  \emph{Combinatorics
of train tracks}, Annals of Mathematics Studies 125, Princeton
University Press

\bibitem{R} \textbf{R Roberts}, \emph{Taut foliations in punctured
surface bundles, III}, preprint

\bibitem{Th} \textbf{W\,P Thurston},  \emph{Three-dimensional
manifolds, Kleinian groups and hyperbolic geometry},
Bull. Amer. Math. Soc. 6 (1982) 357--381

\bibitem{Wa} \textbf{F Waldhausen},  \emph{On Irreducible 3--manifolds
which are Sufficiently Large}, Ann. of Math.  87 (1968) 56--88

\bibitem{Wu} \textbf{Ying-Qing Wu}, \emph{Sutured manifold
hierarchies, essential laminations, and Dehn surgery}, J. Differential
Geometry, 48 (1998) 407--437

\end{thebibliography}
\end{document}